\documentclass[DIV=14,letterpaper]{scrartcl}

\date{}

\usepackage{tikz}
\usetikzlibrary{calc}
\usetikzlibrary{matrix}
\usepackage{amsmath,amssymb,amsfonts,amsthm}
\usepackage[pdftex]{hyperref}
\usepackage{pdfpages}
\usepackage{graphicx}
\usepackage{subfigure}
\allowdisplaybreaks
\usepackage{algorithm,algpseudocode,float}
\usepackage{mathabx}

\newcommand{\bey}{\begin{eqnarray}}
\newcommand{\eey}{\end{eqnarray}}

\newcommand{\beq}{\begin{equation}}
\newcommand{\eeq}{\end{equation}}
\theoremstyle{plain}
\newtheorem{thm}{\hspace{6mm}Theorem}[section]

\theoremstyle{definition}

\newtheorem{exam}{\hspace{6mm}Example}[section]

\title{Meshfree finite difference solution of homogeneous Dirichlet problems of the fractional Laplacian}

\author{
Jinye Shen\thanks{School of Mathematics, Southwestern University of Finance and Economics,
Chengdu, Sichuan, China. {\em jyshen@swufe.edu.cn}},
\;
Bowen Shi\thanks{School of Mathematics, Southwestern University of Finance and Economics,
Chengdu, Sichuan, China. {\em BW001209@163.com}},
\; and\,
Weizhang Huang\thanks{Department of Mathematics, University of Kansas, Lawrence, Kansas, U.S.A.
{\em whuang@ku.edu}}
}

\begin{document}
\vskip 1cm
\maketitle

\begin{abstract}
A so-called grid-overlay finite difference method (GoFD) was proposed recently for the numerical solution of homogeneous Dirichlet boundary value problems of the fractional Laplacian on arbitrary bounded domains. It was shown to have advantages of both finite difference and finite element methods,
including its efficient implementation through the fast Fourier transform and ability to work for complex domains and with mesh adaptation. The purpose of this work is to study GoFD in a meshfree setting, a key to which is to construct the data transfer matrix from a given point cloud to a uniform grid. Two approaches are proposed, one based on the moving least squares fitting and the other based on the Delaunay triangulation and piecewise linear interpolation. Numerical results obtained for examples with convex and concave domains and various types of point clouds are presented. They show that both approaches lead to comparable results. Moreover, the resulting meshfree GoFD converges at a similar order as GoFD with unstructured meshes and finite element approximation as the number of points in the cloud increases. Furthermore, numerical results show that the method is robust to random perturbations in the location of the points.
\end{abstract}

\noindent
\textbf{AMS 2020 Mathematics Subject Classification.} 65N06, 35R11

\noindent
\textbf{Key Words.} Fractional Laplacian, meshfree, finite difference, arbitrary domain, overlay grid

\newpage

\section{Introduction}

The fractional Laplacian is a fundamental non-local operator in the modeling of anomalous dynamics and its numerical
approximation has attracted considerable attention recently; e.g. see \cite{Antil-2022,Huang2016,Lischke-2020} and references therein.
A number of numerical methods have been developed, for example,
finite difference (FD) methods
\cite{DuNing2019,Duo-2018,Hao2021,Huang2014,Huang2016,Ying-2020,Pang-2012,Sunjing2021,Wang2012},
finite element methods
\cite{Acosta2017,Acosta201701,Ainsworth-2017,Ainsworth-2018,Bonito2019,Faustmann2022,Tian2013},
spectral methods \cite{Song2017,LiHuiyuan2022},
discontinuous Galerkin methods \cite{Du2019},
meshfree methods \cite{Burkardt-2021,Pang-2015}, and
sinc-based methods \cite{Antil-2021}.
Loosely speaking, methods such as FD methods are constructed on uniform grids and have the advantage
of efficient matrix-vector multiplication via the fast Fourier transform (FFT) but do not work for complex domains
and with mesh adaptation.
On the other hand, methods such as finite element methods
can work for arbitrary bounded domains and with mesh adaptation
but suffer from slowness of stiffness matrix assembling and matrix-vector multiplication because the stiffness matrix is a full matrix.
Few methods can do both.
A sparse approximation to the stiffness matrix of finite element methods and an efficient multigrid implementation
have been proposed by Ainsworth and Glusa \cite{Ainsworth-2017,Ainsworth-2018}.

\begin{figure}[ht!]
\centering
\includegraphics[width=0.22\textwidth]{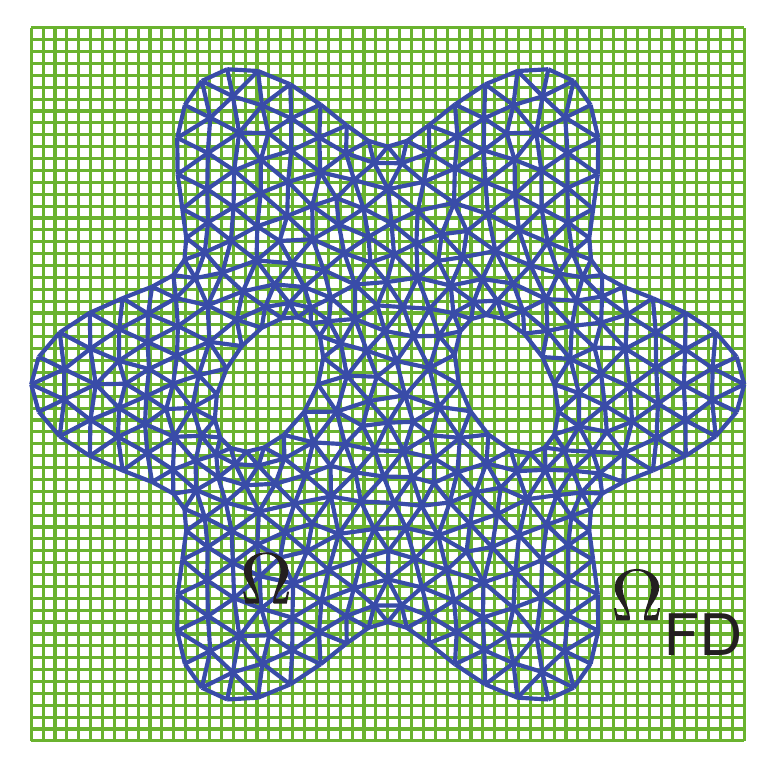}
\caption{A rectangular/cubic domain $\Omega_{\text{FD}} \supseteq \Omega$ and a uniform grid $\mathcal{T}_{\text{FD}}$ (in green color)
on $\Omega_{\text{FD}}$ are used in the GoFD solution of fractional Laplacian BVPs. BVPs are solved on an unstructured mesh $\mathcal{T}_h$ (in blue color).}
\label{fig:gridoverlay-1}
\end{figure}

A so-called grid-overlay finite difference method (GoFD) was recently proposed in \cite{Huang2023}
for the numerical solution of homogeneous Dirichlet boundary value problems (BVPs) of the fractional Laplacian.
It combines the advantages of both finite difference and finite element methods to have efficient implementation via FFT
while being able to work for complex domains and with mesh adaptation.
GoFD uses an unstructured mesh for the underlying domain and a uniform grid overlaying it (cf. Fig.~\ref{fig:gridoverlay-1})
and then constructs an FD approximation for the fractional Laplacian by combining a uniform-grid FD approximation
of the fractional Laplacian and a data transfer between the uniform grid and the unstructured mesh.
A sufficient condition was provided for the solvability of
the GoFD approximation and the convergence order of GoFD
was demonstrated numerically in \cite{Huang2023}.

The objective of this work is to study GoFD in a meshfree/meshless setting in two dimensions. Meshfree methods
offer a more convenient approach to handle complex geometries than finite difference and finite element methods
by using point clouds instead of meshes; e.g., see \cite{Liu-GR-2010,Suchde-2023,Trobec2015}.
Recall that a mesh consists of a set of points and the connectivity between them.
Generating a set of points (called a point cloud) is much easier than generating a mesh for a domain with complex geometry.
Here, we explore the feasibility of using GoFD with point clouds for the numerical solution
of homogeneous Dirichlet problems of the fractional Laplacian. In the original development of GoFD in \cite{Huang2023},
an unstructured mesh is used for the underlying domain although its connectivity is used only in the construction
of the transfer matrix associated with element-wise linear interpolation.
A key to the success of GoFD in a meshfree setting is then to construct the transfer matrix for a given point cloud.
In this work we consider two approaches, one based on the moving least squares fitting
with inverse distance weighting and the other based on the Delaunay triangulation and piecewise linear interpolation.
The former is completely meshfree while the latter is not. Nevertheless, we still consider
the Delaunay triangulation approach as a worthwhile option because the Delaunay triangulation has been widely and successfully
used for scattered data interpolation in data processing. Moreover, constructing the transfer matrix for piecewise linear interpolation
on a mesh does not require the mesh to have a high quality as in finite difference and finite element simulation where
the conditioning of the stiffness matrix and/or the approximation accuracy can deteriorate significantly when
the quality of the mesh worsens. The Delaunay triangulation (without subsequent post-processing or optimization)
provides an efficient way to generate a crude mesh for constructing the transfer matrix for GoFD.
Numerical examples with convex and concave domains and various types of point clouds
will be presented to demonstrate  the accuracy and robustness of the resulting meshfree
GoFD with both approaches to construct the transfer matrix.

It is worth pointing out that GoFD uses a uniform grid (cf. Fig.~\ref{fig:gridoverlay-1}). The generation of this grid
has nothing to do with the domain boundary and is almost cost-free. Thus, this grid cannot be compared
with background grids used in several meshfree methods for the computation of integrals needed in weak formulations \cite{Suchde-2023}.

An online of the paper is as follows. GoFD with unstructured meshes for $\Omega$ will be described in Section~\ref{SEC:GoFG}.
The approaches based on the moving least squares fitting and Delaunay triangulation to construct the transfer matrix for a given point
cloud will be discussed in Section~\ref{SEC:IhFD}. Numerical examples with convex and concave domains and various types of
point clouds are presented in Section \ref{SEC:numerics}. Finally, Section~\ref{SEC:conclusions}
contains conclusions and further remarks.

\section{The grid-overlay finite difference method}
\label{SEC:GoFG}

We consider the homogeneous Dirichlet BVP
\begin{equation}
\label{BVP-1}
\begin{cases}
(-\Delta)^{s}  u = f, & \quad \text{ in } \Omega ,
\\
u = 0, & \quad \text{ in } \Omega^c ,
\end{cases}
\end{equation}
where $\Omega$ is a bounded domain in two dimensions,
$ \Omega^c \equiv \mathbb{R}^2\setminus \Omega$ is the complement of $\Omega$, $f$ is a given function, and
$(-\Delta)^{s}$ is the fractional Laplacian of order $s \in (0,1)$. The fractional Laplacian can be expressed in terms of singular integrals as
\begin{equation}
\label{FL-1}
(-\Delta)^s u(\mbox{\boldmath \( x \)})
=C_{2,s} \text{ p.v. } \int_{\mathbb{R}^2}
\frac{u(\mbox{\boldmath \( x \)})-u(\mbox{\boldmath \( y \)})}{ \mid \mbox{\boldmath \( x \)}-\mbox{\boldmath \( y \)} \mid ^{2+2s}}d \mbox{\boldmath \( y \)},\quad C_{2,s}=\frac{2^{2s}s \Gamma(s+1)}{\pi \Gamma (1-s)}
\end{equation}
or in terms of the Fourier transform as
\begin{equation}
\label{FL-2}
(-\Delta)^s u = \mathcal{F}^{-1}( \mid \mbox{\boldmath \( \xi \)} \mid^{2s} \mathcal{F}(u)),
\end{equation}
where p.v. stands for the Cauchy principal value, $\Gamma(\cdot)$ is the gamma function, and  $\mathcal{F}$ and $\mathcal{F}^{-1}$
denote the Fourier transform and the inverse Fourier transform, respectively.

\subsection{Uniform-grid FD approximation of the fractional Laplacian}
\label{SEC:FD}

We start the description of GoFD for solving (\ref{BVP-1}) with an FD approximation of the fractional Laplacian on a uniform grid.
Uniform-grid FD approximation of the fractional Laplacian has been studied extensively; for example, see
\cite{Hao2021,Huang2023,Huang2016,Ortigueira2006,Ortigueira2008}.
For our current situation, we first choose a rectangular domain $\Omega_{\text{FD}}$ with half width $R_{\text{FD}}$ such that
it covers $\Omega$; see Fig.~\ref{fig:gridoverlay-1}.
Then, we define a uniform grid $\mathcal{T}_{\text{FD}}$ with spacing $h_{\text{FD}}$ as
\[
\mathcal{T}_{\text{FD}}: \quad \quad (x_j, y_k) \equiv (j\, h_{\text{FD}}, k\, h_{\text{FD}}), \quad j, k = 0, \pm 1, \pm 2, ... .
\]
The choice of $h_{\text{FD}}$ will be discussed in the next subsection (cf. Theorem~\ref{thm:Ah-1}).
For the moment, we consider $\mathcal{T}_{\text{FD}}$ as an infinite uniform grid. But we will use
its restriction on $\Omega_{\text{FD}}$ in actual computation (cf. (\ref{FDFL-2}) below).

We consider the uniform-grid FD approximation of the fractional Laplacian
based on its Fourier transform representation (\ref{FL-2}). In two dimensions, (\ref{FL-2}) reduces to
\[
(-\Delta)^s u (x,y) = \frac{1}{(2\pi)^2} \iint_{\mathbb{R}^2} (\xi^2+\eta^2)^s \hat{u}(\xi,\eta) e^{i x \xi + i y \eta} d \xi d \eta,
\]
where
\[
\hat{u}(\xi,\eta) = \iint_{\mathbb{R}^2} u(x,y) e^{- i x \xi - i y \eta} d x d y .
\]
Define the discrete Fourier transform of $u$ on $\mathcal{T}_{\text{FD}}$ as
\[
\check{u}(\xi,\eta) =
\sum_{m=-\infty}^{\infty}\sum_{n=-\infty}^{\infty}
u_{m,n} e^{i x_m \xi+i y_n \eta},
\]
where $u_{m,n} = u(x_m,y_n)$.
Applying the discrete Fourier transform to the five-point central finite difference approximation of the Laplacian,
we can obtain the discrete Fourier transform of the FD approximation of the fractional Laplacian as
\[
\check{(-\Delta_h)^s u} (\xi, \eta)  = \frac{1}{h_{\text{FD}}^{2 s}}
\left (4 \sin^{2}(\frac{\xi h_{\text{FD}}}{2}) + 4 \sin^{2}(\frac{\eta h_{\text{FD}}}{2})\right )^{s} \check{u}(\xi,\eta) .
\]
Then, the FD approximation of the fractional Laplacian is given by
\begin{align*}
(-\Delta_h)^s u(x_j,y_k) & = \frac{h_{\text{FD}}^2}{ (2\pi)^2}
\int_{-\frac{\pi}{h_{\text{FD}}}}^{\frac{\pi}{h_{\text{FD}}}} \int_{-\frac{\pi}{h_{\text{FD}}}}^{\frac{\pi}{h_{\text{FD}}}}
\check{(-\Delta_h)^s u} (\xi,\eta) e^{i x_j \xi+i y_k \eta} d \xi d\eta .
\end{align*}
We can rewrite this into
\begin{align}
(-\Delta_h)^s u(x_j,y_k)  = \frac{1}{h_{\text{FD}}^{2 s}} \sum_{m=-\infty}^{\infty} \sum_{n=-\infty}^{\infty} A_{(j,k),(m,n)} u_{m,n} ,
\label{FDFL-1}
\end{align}
where
\begin{align}
\label{A-1}
A_{(j,k),(m,n)} = T_{j-m,k-n},
\end{align}
and
\begin{align}
\label{T-1}
T_{p,q}  = \frac{1}{(2\pi)^2} \int_{-\pi}^{\pi} \int_{-\pi}^{\pi}
\left [ 4 \sin^2\left (\frac{\xi}{2}\right ) + 4 \sin^2\left (\frac{\eta}{2}\right )\right ]^s e^{i p \xi+i q \eta} d \xi d \eta .
\end{align}
Notice that $T_{p,q}$'s are the coefficients of the Fourier series expansion of
$(4 \sin^2({\xi}/{2}) + 4 \sin^2({\eta}/{2}))^s$. Moreover, it is not difficult to show
\begin{equation}
\label{T-3}
T_{-p,-q} = T_{p,q}, \quad T_{-p,q} = T_{p,q}, \quad T_{p,-q} = T_{p,q} .
\end{equation}

 For the numerical solution of (\ref{BVP-1}),  we focus on functions vanishing in $\Omega^c$ (and therefore in $\Omega_{\text{FD}}^c$).
 For those functions, (\ref{FDFL-1}) reduces to a finite summation, i.e.,
\begin{align}
(-\Delta_h)^s u(x_j,y_k) =
\frac{1}{h_{\text{FD}}^{2 s}} \sum_{m=-N_{\text{FD}}}^{N_{\text{FD}}} \sum_{n=-N_{\text{FD}}}^{N_{\text{FD}}} A_{(j,k),(m,n)} u_{m,n} ,
\; -N_{\text{FD}} \le j, k \le N_{\text{FD}} ,
\label{FDFL-2}
\end{align}
where $N_{\text{FD}} = [R_{\text{FD}}/h_{\text{FD}}] + 1$, where $[R_{\text{FD}}/h_{\text{FD}}]$ denotes the integer part of
$R_{\text{FD}}/h_{\text{FD}}$.
Hereafter, this FD approximation of the fractional Laplacian (after omitting the factor $1/h_{\text{FD}}^{2 s}$) will be referred to as
$A_{\text{FD}}$. It is known \cite{Huang2023} that $A_{\text{FD}}$ is symmetric and positive definite.

$T_{p,q}$'s can be computed using the rectangle quadrature rule and FFT (e.g., see \cite{Huang2023}). Moreover, we have
\begin{align*}
& \sum_{m=-N_{\text{FD}}}^{N_{\text{FD}}} \sum_{n=-N_{\text{FD}}}^{N_{\text{FD}}} A_{(j,k),(m,n)} u_{m,n}
= \sum_{m = -N_{\text{FD}}}^{N_{\text{FD}}} \sum_{n=-N_{\text{FD}}}^{N_{\text{FD}}} T_{j-m,k-n} u_{m,n}
\\
& = \frac{1}{(4N_{\text{FD}})^2} \sum_{p=-2N_{\text{FD}}}^{2N_{\text{FD}}-1} \sum_{q=-2N_{\text{FD}}}^{2N_{\text{FD}}-1} \check{T}_{p,q} \check{u}_{p,q}
(-1)^{(p+2N_{\text{FD}}) + (q+2N_{\text{FD}})}
\\
& \qquad \qquad  \qquad \qquad \qquad \qquad
\times e^{\frac{i 2 \pi (p+2N_{\text{FD}})(j+N_{\text{FD}})}{4N_{\text{FD}}}+\frac{i 2 \pi (q+2N_{\text{FD}})(k+N_{\text{FD}})}{4N_{\text{FD}}}} ,
\end{align*}
where
\begin{align*}
\check{T}_{p,q} = \sum_{m=-2N_{\text{FD}}}^{2N_{\text{FD}}-1} \sum_{n=-2N_{\text{FD}}}^{2N_{\text{FD}}-1} T_{m,n} e^{-\frac{i 2 \pi (m+2N_{\text{FD}})(p+2N_{\text{FD}})}{4N_{\text{FD}}}
-\frac{i 2 \pi (n+2N_{\text{FD}})(q+2N_{\text{FD}})}{4N_{\text{FD}}}},\\
 - 2N_{\text{FD}} \le p, q \le 2 N_{\text{FD}} ,
\\
\check{u}_{p,q}  = \sum_{m = -N_{\text{FD}}}^{N_{\text{FD}}} \sum_{n=-N_{\text{FD}}}^{N_{\text{FD}}} u_{m,n} e^{-\frac{i 2 \pi (m+2N_{\text{FD}})(p+2N_{\text{FD}})}{4N_{\text{FD}}}
-\frac{ i 2 \pi (n+2N_{\text{FD}})(q+2N_{\text{FD}})}{4N_{\text{FD}}}}, \\
 - 2N_{\text{FD}} \le p, q \le 2 N_{\text{FD}} .
\end{align*}
Thus, $\check{T}_{p,q}$'s, $\check{u}_{p,q}$'s, and the multiplication of $A_{\text{FD}}$
with vectors can be computed using FFT.

\subsection{GoFD for homogeneous Dirichlet BVPs}
\label{SEC:GoFD-0}

Now we describe GoFD for BVP (\ref{BVP-1}).
To start with, we assume that an unstructured simplicial mesh $\mathcal{T}_h$ is given for $\Omega$.
(The mesh $\mathcal{T}_h$ will be replaced by a set of points on $\Omega$ in the meshfree setting; see the next section.)
Denote the vertices of $\mathcal{T}_h$ by $(x_j,y_j)$, $j = 1, ..., N_{vi}, ..., N_v$, where $N_{v}$
and $N_{vi}$ are the numbers of the vertices and interior vertices,  respectively. Here, we assume that
the vertices are arranged so that the interior vertices are listed before the boundary ones.
Next, we take $\Omega_{\text{FD}} \equiv (-R_{\text{FD}},R_{\text{FD}})^d$ to cover $\Omega$
and create an overlaying uniform grid (denoted by $\mathcal{T}_{\text{FD}}$)
with $2 N_{\text{FD}}+1$ nodes in each axial direction for some positive integer $N_{\text{FD}}$ (cf. Fig.~\ref{fig:gridoverlay-1}).
The vertices of $\mathcal{T}_{\text{FD}}$ are denoted by $(x_k^{\text{FD}},y_k^{\text{FD}})$, $k = 1, ..., N_{v}^{\text{FD}}$.
Then, a uniform-grid FD approximation $h_{\text{FD}}^{-2 s} A_{\text{FD}}$ can be constructed
for the fractional Laplacian on $\mathcal{T}_{\text{FD}}$ as described in the previous subsection.
The GoFD approximation for BVP \eqref{BVP-1} is then defined as
\begin{equation}
\label{GoFD-1}
\frac{1}{h_{\text{FD}}^{2 s}} D_h^{-1} (I_{h}^{\text{FD}})^T A_{\text{FD}} I_{h}^{\text{FD}} \mbox{\boldmath \( u \)}_h
=  \mbox{\boldmath \( f \)}_h,
\end{equation}
where $\mbox{\boldmath \( u \)}_h = ( u_h(x_1,y_1), ..., u_h(x_{N_{vi}},y_{N_{vi}}))^T$, $\mbox{\boldmath \( f \)}_{h}=( f(x_1,y_1), ..., f(x_{N_{vi}},y_{N_{vi}}))^T$,
$I_h^{\text{FD}}$ is a transfer matrix from $\mathcal{T}_h$ to $\mathcal{T}_{\text{FD}}$,
and $D_h$ is the diagonal matrix formed by the column sums of $I_h^{\text{FD}}$.
The invertibility of $D_h$ will be addressed in Theorem~\ref{thm:Ah-1} below.

We consider the transfer matrix for piecewise linear interpolation.
For any function $u(x,y)$, we can express its piecewise linear interpolation as
\begin{equation}
\label{interp-1}
I_h u (x,y) = \sum_{j=1}^{N_{v}} u(x_j,y_j) \phi_j(x,y),
\end{equation}
where $\phi_j$ is the Lagrange-type linear basis function associated with vertex $(x_j,y_j)$
and is extended to $\Omega^c$ using $\phi_j(x,y) = 0$ for all $(x,y) \in \Omega^c$.
Then,
\[
I_h u (x_k^{\text{FD}},y_k^{\text{FD}}) = \sum_{j=1}^{N_v} u(x_j,y_j) \phi_j(x_k^{\text{FD}},y_k^{\text{FD}}) , \quad k = 1, ..., N_v^{\text{FD}}.
\]
From this we obtain the entries of $I_h^{\text{FD}}$ as
\begin{equation}
\label{IhFD-1}
(I_h^{\text{FD}})_{k,j} = \phi_j(x_k^{\text{FD}},y_k^{\text{FD}}), \quad k = 1, ..., N_v^{\text{FD}}, \; j = 1, ..., N_v .
\end{equation}
Notice that $(I_h^{\text{FD}})_{k,j}  = 0$ when $(x_k^{\text{FD}},y_k^{\text{FD}}) \notin \Omega$.

\begin{thm}[\cite{Huang2023}]
\label{thm:Ah-1}
If $N_{\text{FD}}$ is chosen sufficiently large such that
\begin{equation}
h_{\text{FD}} \equiv \frac{R_{\text{FD}}}{N_{\text{FD}}} \le \frac{a_h}{(d+1)\sqrt{d}},
\label{hFD-1}
\end{equation}
where $a_h$ is the minimum element height of $\mathcal{T}_h$, then the transfer matrix $I_h^{\text{FD}}$ associated with piecewise
linear interpolation has full column rank, $D_h$ is invertible, and the GoFD approximation matrix (\ref{GoFD-1})
for the fractional Laplacian is similar to a symmetric and positive definite matrix and thus invertible.
\end{thm}

The above theorem gives a sufficient condition (\ref{hFD-1}) for the invertibility of $A_h$. This condition is on the conservative side
in general. Numerical experiment suggests that a less restrictive condition such as
$ h_{\text{FD}} \le a_h$ can be used.

We now turn our attention to the solution of the linear system (\ref{GoFD-1}). It can be rewritten into a symmetric system as
\begin{equation}
\label{GoFD-2}
(I_{h}^{\text{FD}})^T A_{\text{FD}} I_{h}^{\text{FD}}  \mbox{\boldmath \( u \)}_h =  h_{\text{FD}}^{2 s} D_h \mbox{\boldmath \( f \)}_h .
\end{equation}
Notice that the coefficient matrix of this system is a full matrix and thus direct solution is unrealistic. On the other hand,
its multiplication with vectors can be implemented efficiently. To see this, we recall that
$I_{h}^{\text{FD}}$ is a sparse matrix of size $(2N_{\text{FD}}+1)^2 \times N_v$, where $(2N_{\text{FD}}+1)^2$ and $N_v$
are the numbers of vertices in $\mathcal{T}_{\text{FD}}$ and $\mathcal{T}_h$, respectively. From (\ref{IhFD-1}) and the locality
of linear basis functions we can see that the number of
the non-zero entries of $I_{h}^{\text{FD}}$ is $\mathcal{O}(N_{\text{FD}}^2)$. Thus, the multiplication
of $I_{h}^{\text{FD}}$ or $(I_{h}^{\text{FD}})^T$
with a vector can be carried out in $\mathcal{O}(N_{\text{FD}}^2)$ flops.
Moreover, from the previous subsection we have seen that the multiplication of $A_{\text{FD}}$ with a vector can be carried
using FFT with $\mathcal{O}(N_{\text{FD}}^2 \log (N_{\text{FD}}^2))$ flops.
Thus, the number of the operations needed to carry out the multiplication of the coefficient matrix of (\ref{GoFD-2}) with a vector is
\begin{equation}
\label{complexity-1}
\mathcal{O}(N_{\text{FD}}^2 \log (N_{\text{FD}}^2)) .
\end{equation}
{ Notice this is in the same order as those for the uniform-grid finite difference approximation of the fractional Laplacian.}
If the unstructured mesh $\mathcal{T}_h$ is quasi-uniform, we have $N_v = \mathcal{O}(N)$
($N$ is the number of elements in $\mathcal{T}_h$),
and $a_h = \mathcal{O}(h) = \mathcal{O}(N^{-\frac{1}{2}})$. In this case, if we take $h_{\text{FD}} = a_h$
(cf. (\ref{GoFD-2})) we have $N = \mathcal{O}(N_{\text{FD}}^2)$. Consequently, (\ref{complexity-1}) can be written as
\begin{equation}
\label{complexity-2}
\mathcal{O}(N \log (N)) ,
\end{equation}
which indicates the multiplication of the coefficient matrix of (\ref{GoFD-2}) with a vector can be carried out efficiently.
Moreover, this suggests that the linear system (\ref{GoFD-2}) should be solved with a Krylov subspace method.
The conjugate gradient method (CG) is used in our computation.

{ It is interesting to compare (\ref{complexity-2}) with that of a finite element implementation (e.g., see \cite{Acosta2017}).
Since the stiffness matrix of the finite element approximation to the fractional Laplacian is a full matrix,
the number of the operations needed to carry out its multiplication with a vector is
\begin{equation}
\label{complexity-3}
\mathcal{O}(N_v^2) = \mathcal{O}(N^2),
\end{equation}
which is much more costly than (\ref{complexity-2}).
}

\section{Transfer matrices}
\label{SEC:IhFD}

In this section we study the construction of the transfer matrix $I_h^{\text{FD}}$ in a meshfree setting where only a point cloud
(say, $(x_j, y_j), \; j = 1, ..., N_v$) is assumed to be given on $\bar{\Omega} \equiv \Omega \cup \partial \Omega$.

Recall that an unstructured mesh $\mathcal{T}_h$ for $\Omega$ is used in GoFD described in the previous section
and its connectivity is used only when constructing the transfer matrix $I_h^{\text{FD}}$ (piecewise linear interpolation) from
$\mathcal{T}_h$ to $\mathcal{T}_{\text{FD}}$. Thus, if we can construct $I_h^{\text{FD}}$
from the point cloud $(x_j, y_j), \; j = 1, ..., N_v$ to the set of points $(x_k^{\text{FD}},y_k^{\text{FD}}),\; k = 1, ..., N_v^{\text{FD}}$
without using the connectivity among the points, then GoFD can go with a point cloud on $\bar{\Omega}$ and in this sense, becomes
meshfree or meshless \cite{Trobec2015}. In the following we discuss two approaches to construct $I_h^{\text{FD}}$.

\subsection{Moving least squares fitting with inverse-distance weighting}

In this approach, for any grid point $(x_k^{\text{FD}},y_k^{\text{FD}})$ in $\Omega$, we choose $n$ nearest points
in the set of points $(x_j, y_j)$, where $n$ is a given positive integer. Denote these $n$ points by $(x_{k_j},y_{k_j}),\, j = 1, ..., n$.
Then, a linear polynomial $p(x,y) = a_0^k + a_1^k (x-x_k^{\text{FD}}) + a_2^k (y-y_k^{\text{FD}})$
is obtained using least squares fitting with inverse-distance weighting, i.e.,
\begin{equation}
\label{mls-1}
\min\limits_{a_0^k, a_1^k, a_2^k} \sum_{j=1}^{n}
\left (u(x_{k_j},y_{k_j})-p(x_{k_j},y_{k_j})\right )^2 w_{k_j},
\end{equation}
where $w_{k_j} = 1/\| (x_{k_j},y_{k_j}) - (x_k^{\text{FD}},y_k^{\text{FD}}) \|$ and $\| (x_{k_j},y_{k_j}) - (x_k^{\text{FD}},y_k^{\text{FD}}) \|$
denotes the Euclidean distance between $(x_{k_j},y_{k_j})$ and
$(x_k^{\text{FD}},y_k^{\text{FD}})$. The linear system for $a_0^k$, $a_1^k$, and $a_2^k$ can be obtained as
\begin{align}
& \begin{bmatrix} \sum\limits_{j=1}^n w_{k_j}, & \sum\limits_{j=1}^n w_{k_j} (x_{k_j} -x_k^{\text{FD}}), & \sum\limits_{j=1}^n w_{k_j} (y_{k_j} -y_k^{\text{FD}})\\
\sum\limits_{j=1}^n w_{k_j} (x_{k_j} -x_k^{\text{FD}}), & \sum\limits_{j=1}^n w_{k_j} (x_{k_j} -x_k^{\text{FD}})^2,
& \sum\limits_{j=1}^n w_{k_j} (x_{k_j} -x_k^{\text{FD}})(y_{k_j} -y_k^{\text{FD}}) \\
\sum\limits_{j=1}^n w_{k_j} (y_{k_j} -y_k^{\text{FD}}), & \sum\limits_{j=1}^n w_{k_j} (x_{k_j} -x_k^{\text{FD}}) (y_{k_j} -y_k^{\text{FD}}),
& \sum\limits_{j=1}^n w_{k_j} (y_{k_j} -y_k^{\text{FD}})^2 \end{bmatrix}\notag \\
&\begin{bmatrix} a_0^k \\ a_1^k \\ a_2^k \end{bmatrix}
 =
\begin{bmatrix} \sum\limits_{j=1}^n w_{k_j} u(x_{k_j},y_{k_j}) \\ \sum\limits_{j=1}^n w_{k_j} u(x_{k_j},y_{k_j}) (x_{k_j} -x_k^{\text{FD}})
\\ \sum\limits_{j=1}^n w_{k_j} u(x_{k_j},y_{k_j}) (y_{k_j} -y_k^{\text{FD}}) \end{bmatrix} .
\label{mls-ls-1}
\end{align}
From this, it is clear that $p(x_k^{\text{FD}},y_k^{\text{FD}}) = a_0^k$ and $a_0^k$ is a linear function of $u(x_{k_j},y_{k_j})$, $j = 1, ..., n$.
Denote this relation as $a_0^k = \sum_{j=1}^{n} a_{0,k_j} u_{k_j}$.
Then, we can define the transform matrix through
\[
\sum_{i = 1}^{N_v}  (I_h^{\text{FD}})_{k,i} u(x_i,y_i) = p(x_k^{\text{FD}},y_k^{\text{FD}}) ,
\]
which gives
\begin{equation}
\label{IhFD-2}
(I_h^{\text{FD}})_{k,i} = \begin{cases} a_{0,k_j}, \quad & \text{ for } i = k_j, \; j = 1, ..., n ,\\ 0, \quad &\text{ otherwise}. \end{cases}
\end{equation}
Notice that $(I_h^{\text{FD}})_{k,i} = 0$ when $(x_k^{\text{FD}},y_k^{\text{FD}})$ is not in $\Omega$.

We remark that the inverse-distance weighting in (\ref{mls-1}) gives more weight on the closest points to $(x_k^{\text{FD}},y_k^{\text{FD}})$.
Especially, if a point $(x_{k_j},y_{k_j})$ coincides with $(x_k^{\text{FD}},y_k^{\text{FD}})$, (\ref{mls-1}) gives
$p(x_k^{\text{FD}},y_k^{\text{FD}}) = u_{k_j}$.

In our computation, we choose $n$ nearest points for each grid point $(x_k^{\text{FD}},y_k^{\text{FD}})$
using Matlab's function {\em knnsearch.m} which implements the $k$-nearest neighbors algorithm.

\subsection{Constrained Delaunay triangulation}

In this approach, we propose to use the given point cloud and the boundary information to partition $\Omega$ into
a constrained Delaunay triangulation. A constrained Delaunay triangulation is a generalization of the Delaunay triangulation
forcing certain required segments into the triangulation as edges and can be  computed efficiently \cite{Chew-1989,Shewchuk-2008}.
Once a triangulation is computed, a transform matrix associated with piecewise linear interpolation can be formed
as in Section~\ref{SEC:GoFD-0} (cf.  (\ref{IhFD-1})). Constrained Delaunay triangulation provides an efficient way
to handle concave domains (cf. Examples~\ref{exam-2} and \ref{exam-3})
while (unconstrained) Delaunay triangulation produces a triangulation for the convex hull
of the given point set.

It is worth pointing out that this approach is not completely meshfree.
Nevertheless, it is still important to consider this approach for the numerical solution of (\ref{BVP-1}).
The first reason is that the Delaunay triangulation/tessellation is widely used in data processing to interpolate
from scattered data to an arbitrary set of points. For example, four of the five methods used in Matlab's scattered data interpolation function
{\em griddata.m} are based on data domain triangulation. The second reason is that scattered data interpolation and
in the current situation, construction of the transfer matrix are much less
sensitive to the quality of the triangulation than in finite difference and finite element simulation
where the conditioning of the stiff matrix and/or approximation accuracy can deteriorate significantly when the quality
of the triangulation worsens.

\section{Numerical examples}
\label{SEC:numerics}

In this section we present numerical results obtained for the homogeneous Dirichlet BVP
\begin{equation}
\label{BVP-2}
\begin{cases}
(-\Delta)^s u = 1, \quad & \text{ in } \Omega , \\
u = 0, \quad & \text{ in } \Omega^c ,
\end{cases}
\end{equation}
with various domains and point clouds. We take $h_{\text{FD}}$ as half the minimum distance among the given points
on $\bar{\Omega}$, $n = 5$ in (\ref{mls-1}), and $R_{\text{FD}} = 1.2 R$, where $R$ is the half width of $\Omega$.

\begin{figure}[ht!]
\centering
\subfigure[$s=0.25$]{
\includegraphics[width=0.27\linewidth]{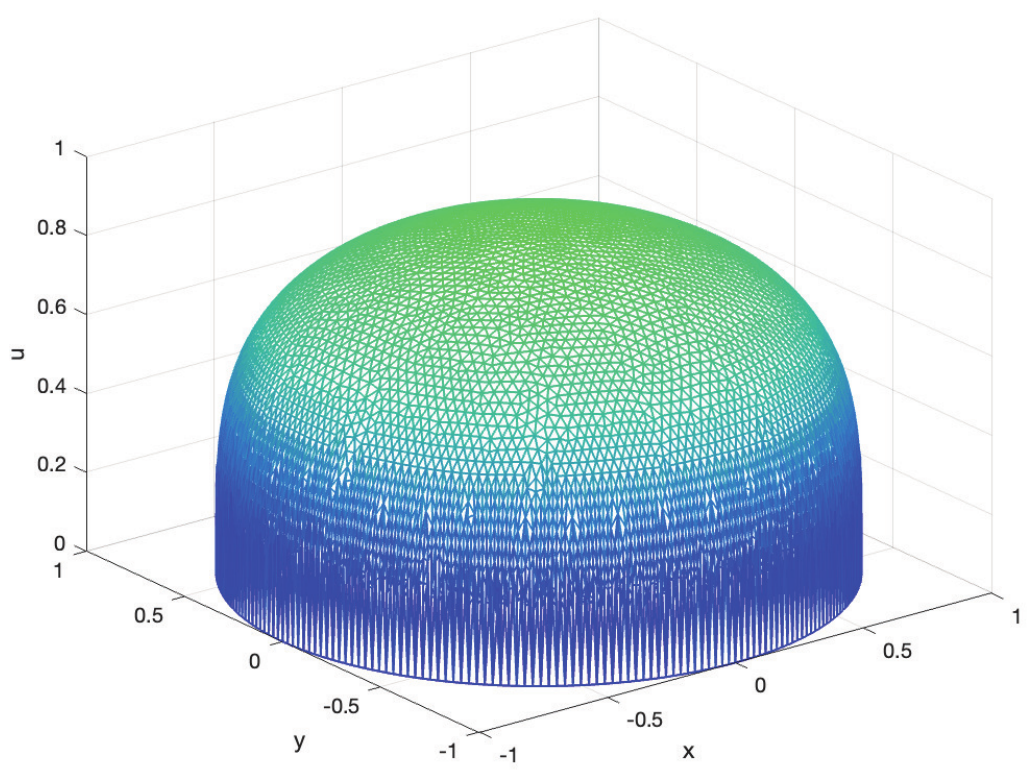}
}
\subfigure[$s=0.5$]{
\includegraphics[width=0.27\linewidth]{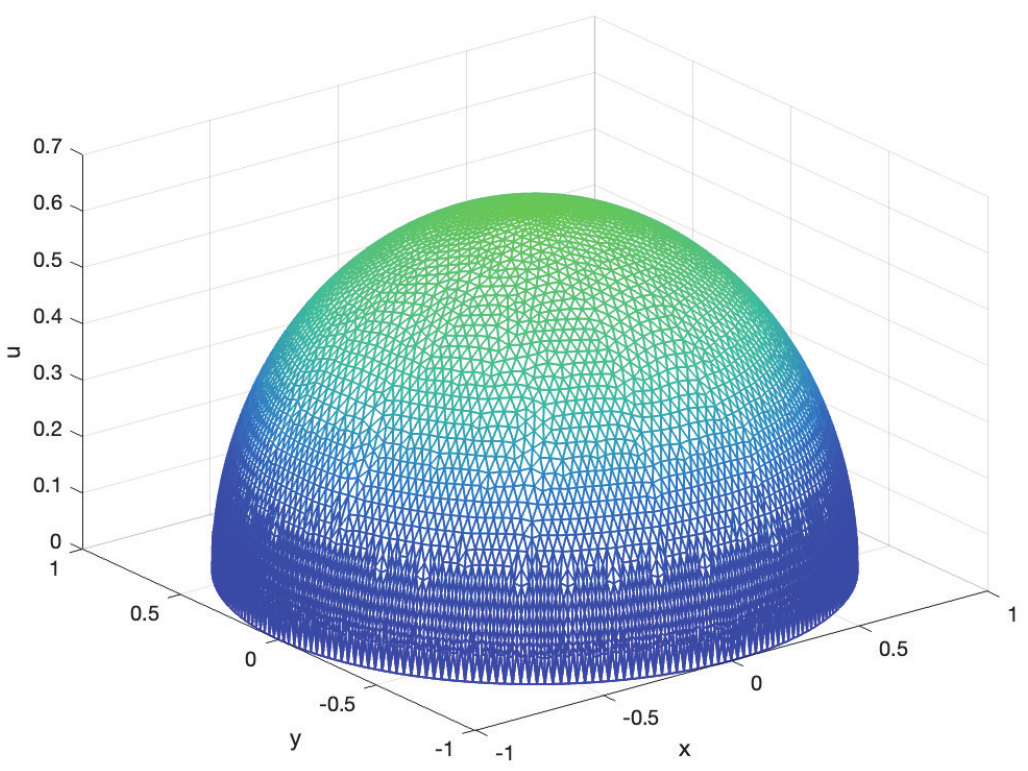}
}
\subfigure[$s=0.75$]{
\includegraphics[width=0.27\linewidth]{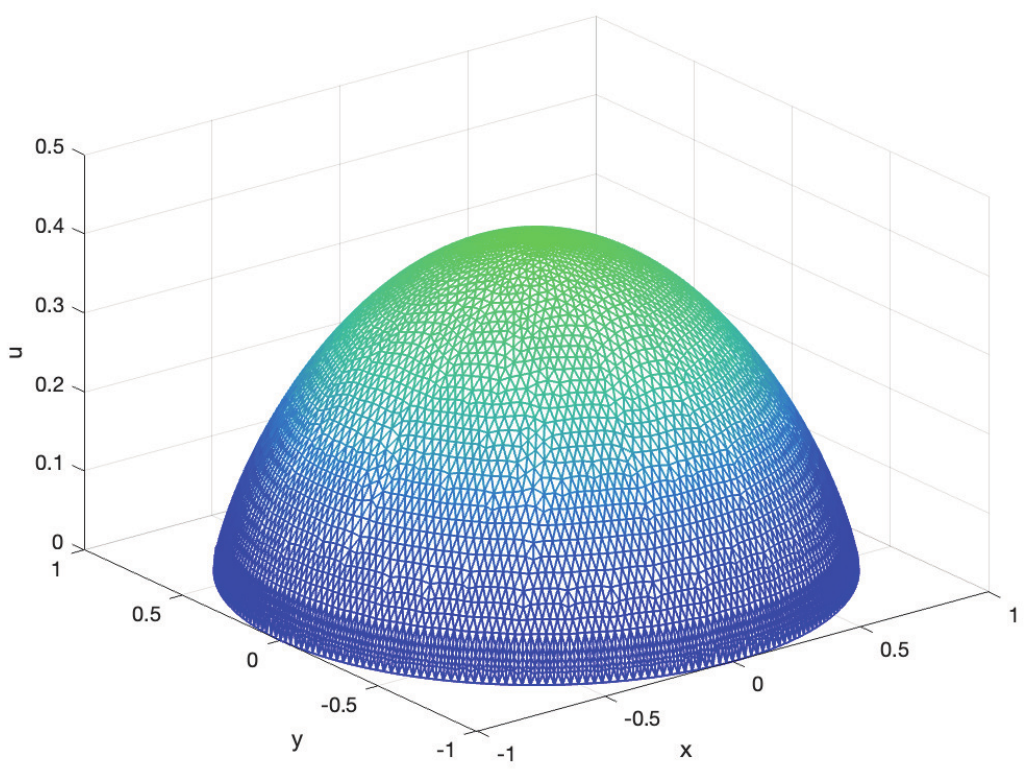}
}
\caption{Example~\ref{exam-1}. The solution (\ref{u-1}) is plotted for $s = 0.25$, 0.5, and 0.75.}
\label{fig:Exam-1-0}
\end{figure}

\begin{figure}[ht!]
\centering
\subfigure[Point Cloud 1]{
\includegraphics[width=0.27\linewidth]{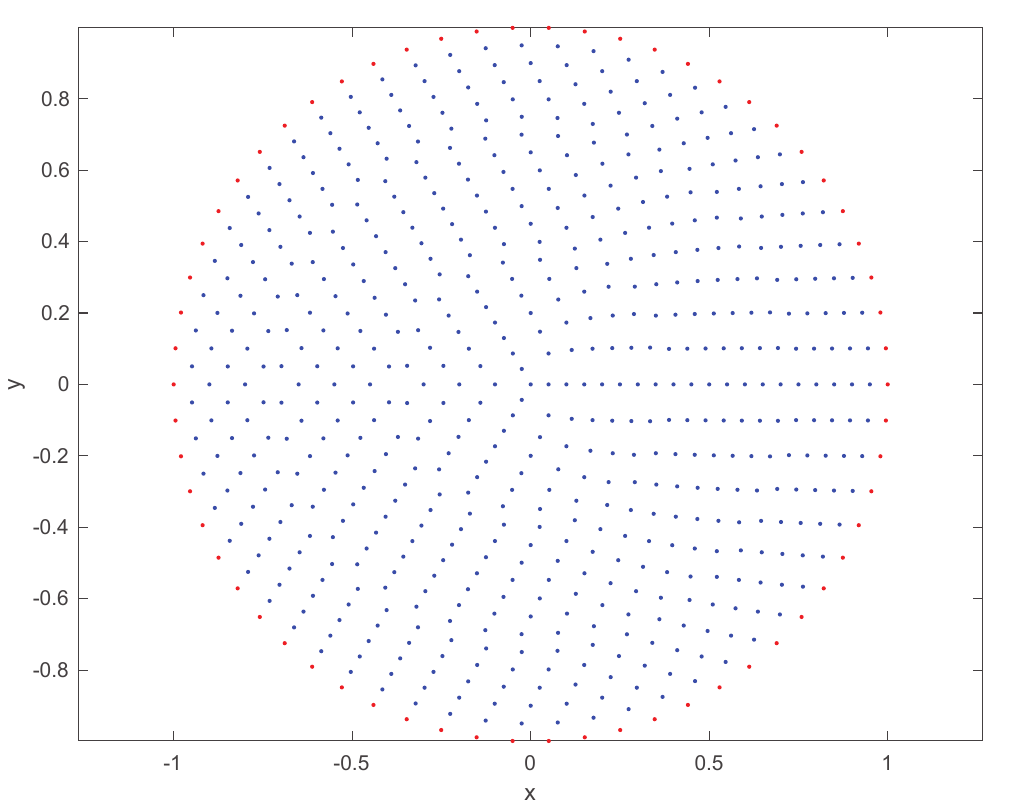}
}
\subfigure[Perturbed Point Cloud 1]{
\includegraphics[width=0.27\linewidth]{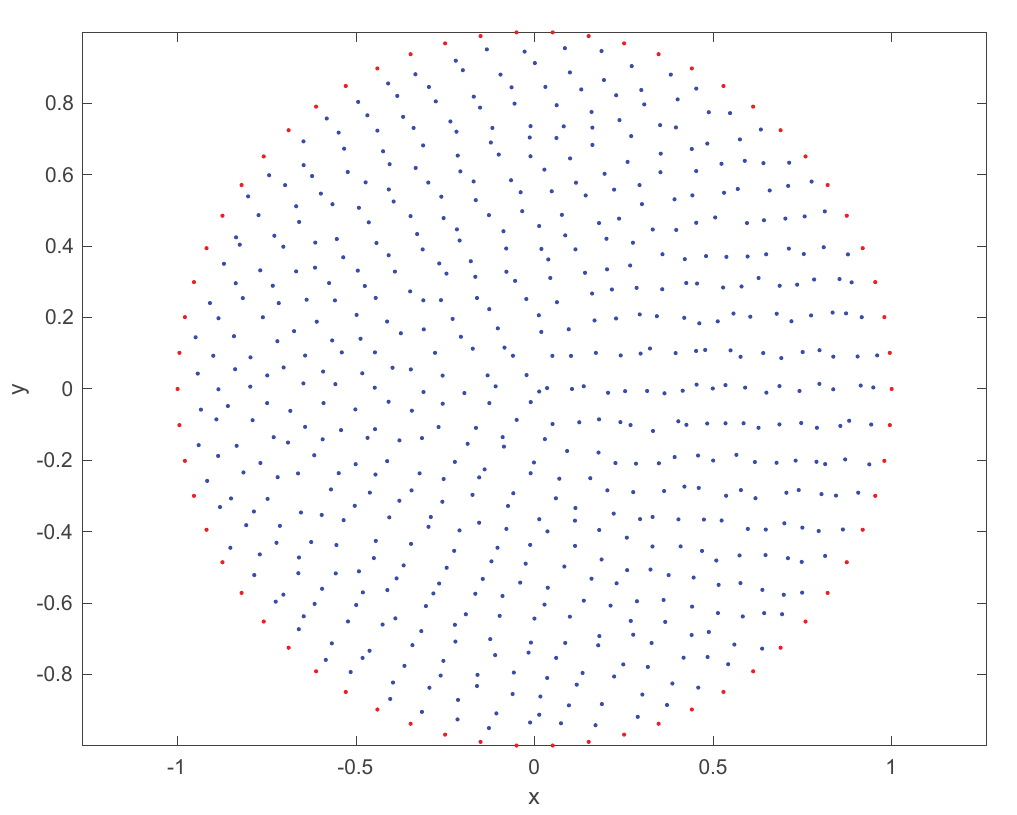}
}
\caption{Example~\ref{exam-1}. Point Cloud 1 and a corresponding perturbed point cloud.}
\label{fig:Exam-1-1}
\end{figure}

\begin{figure}[ht!]
\centering
\subfigure[$s=0.25$]{
\includegraphics[width=0.27\linewidth]{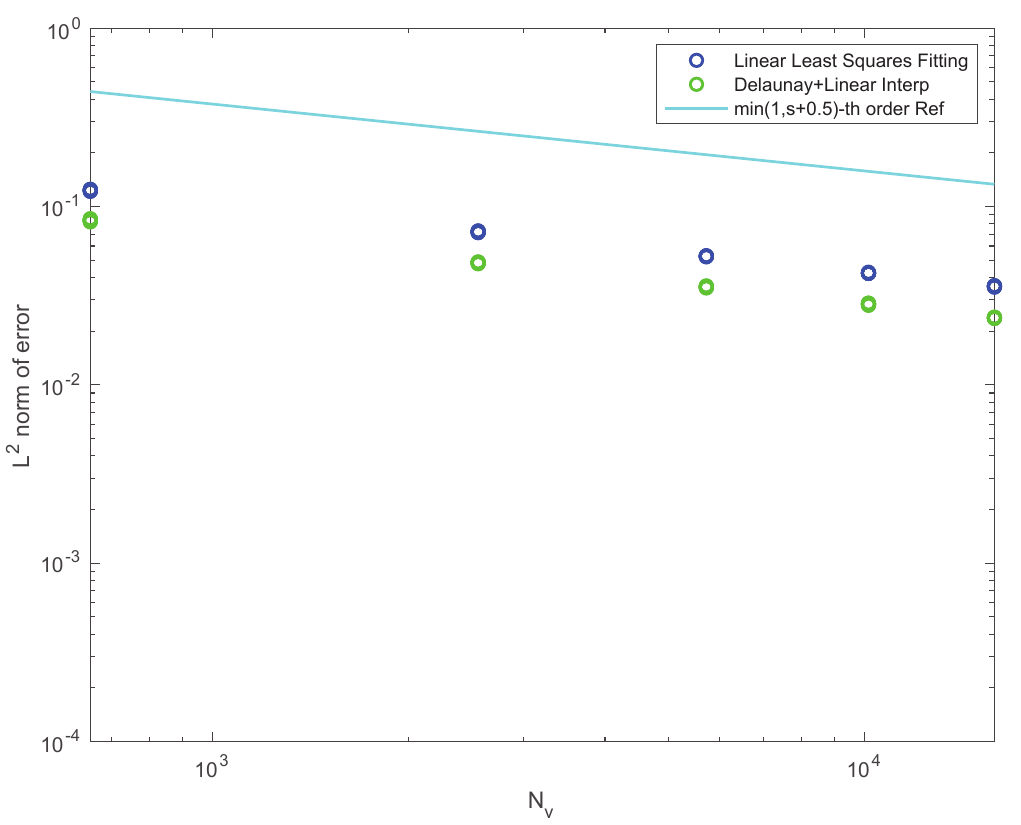}
}
\subfigure[$s=0.5$]{
\includegraphics[width=0.27\linewidth]{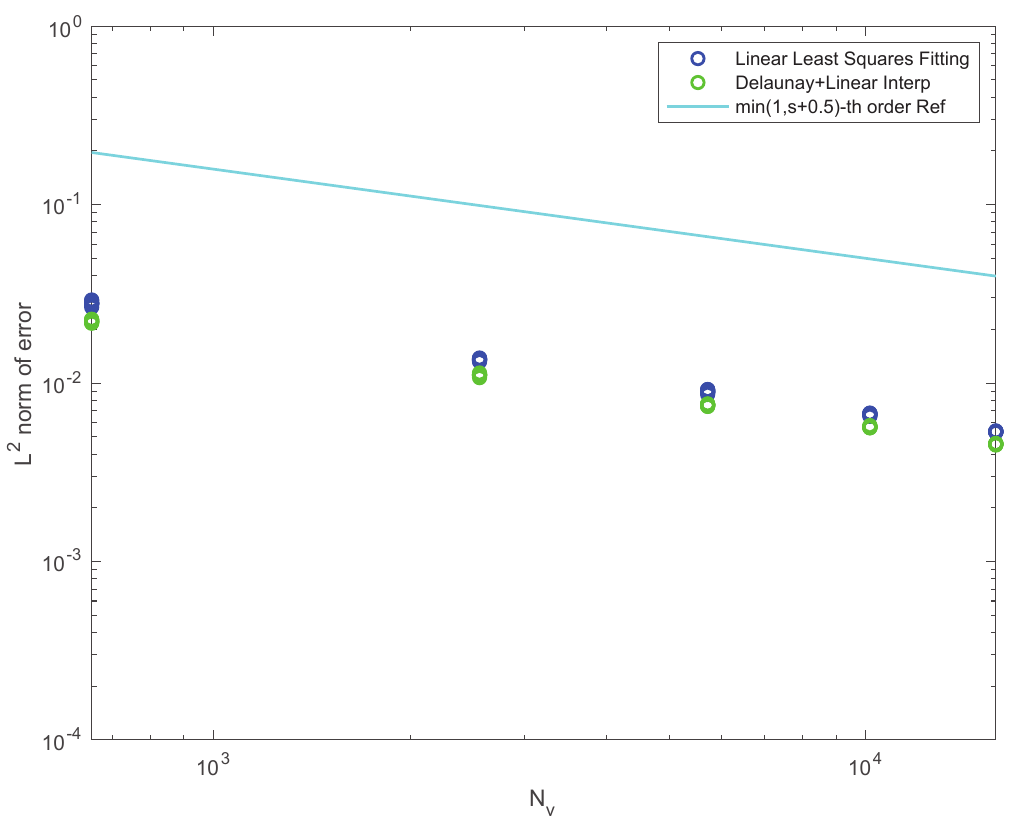}
}
\subfigure[$s=0.75$]{
\includegraphics[width=0.27\linewidth]{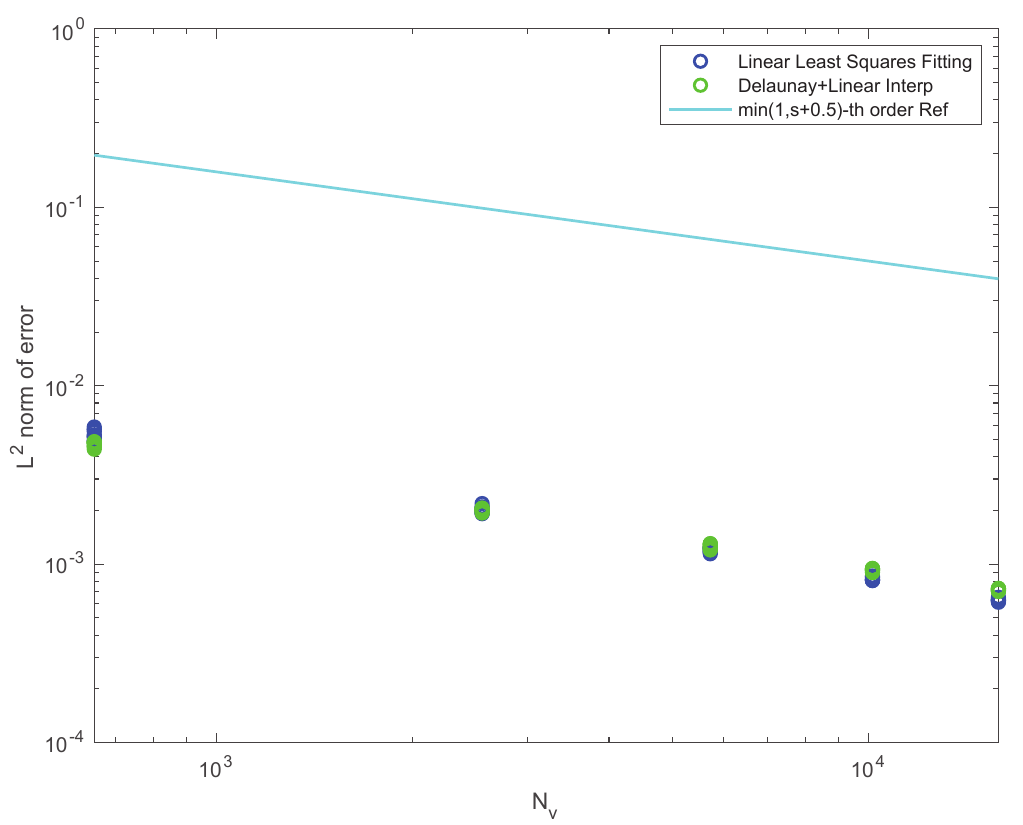}
}
\caption{Example~\ref{exam-1}. The solution error in $L^2$ norm is plotted as a function of $N_v$ for perturbed Point Cloud 1.}
\label{fig:Exam-1-2}
\end{figure}

\begin{figure}[ht!]
\centering
\subfigure[Point Cloud 2]{
\includegraphics[width=0.27\linewidth]{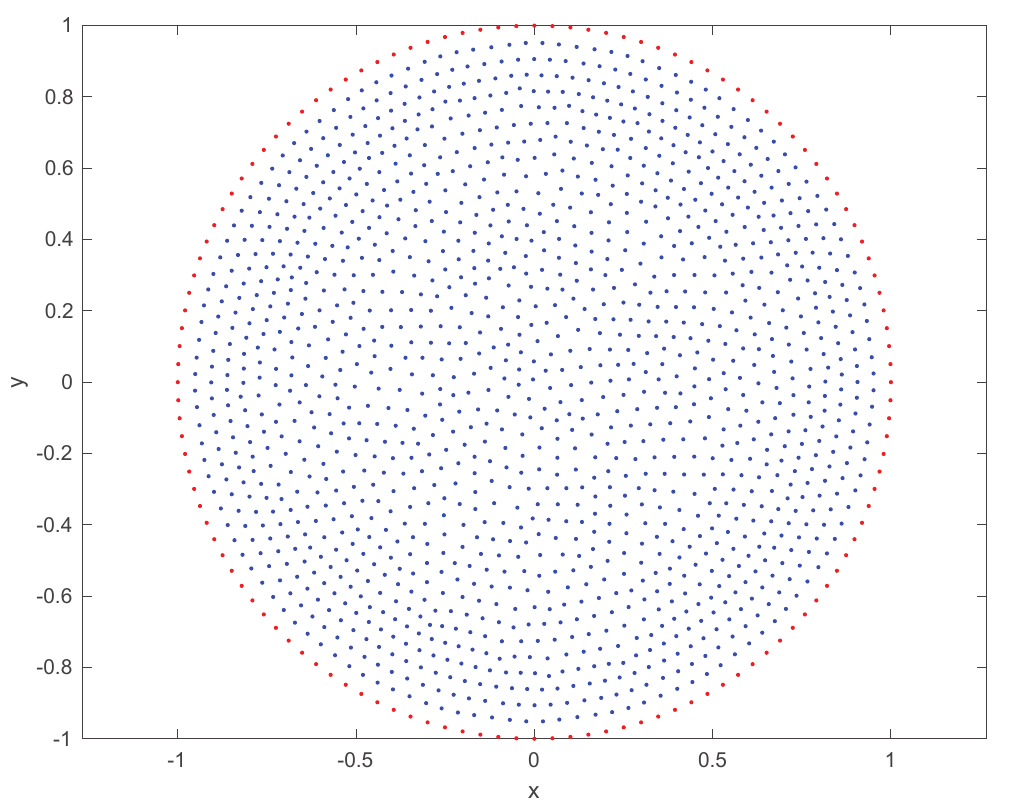}
}
\subfigure[Perturbed Point Cloud 2]{
\includegraphics[width=0.27\linewidth]{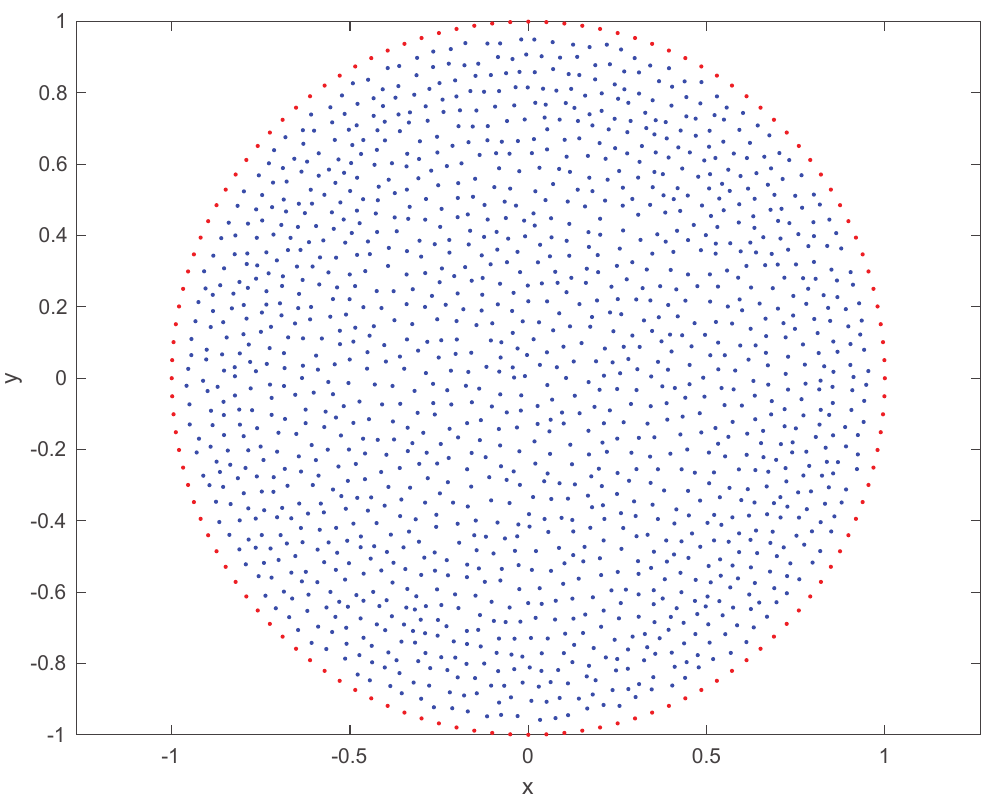}
}
\caption{Example~\ref{exam-1}. Point Cloud 2 and a corresponding perturbed point cloud.}
\label{fig:Exam-1-3}
\end{figure}

\begin{figure}[ht!]
\centering
\subfigure[$s=0.25$]{
\includegraphics[width=0.27\linewidth]{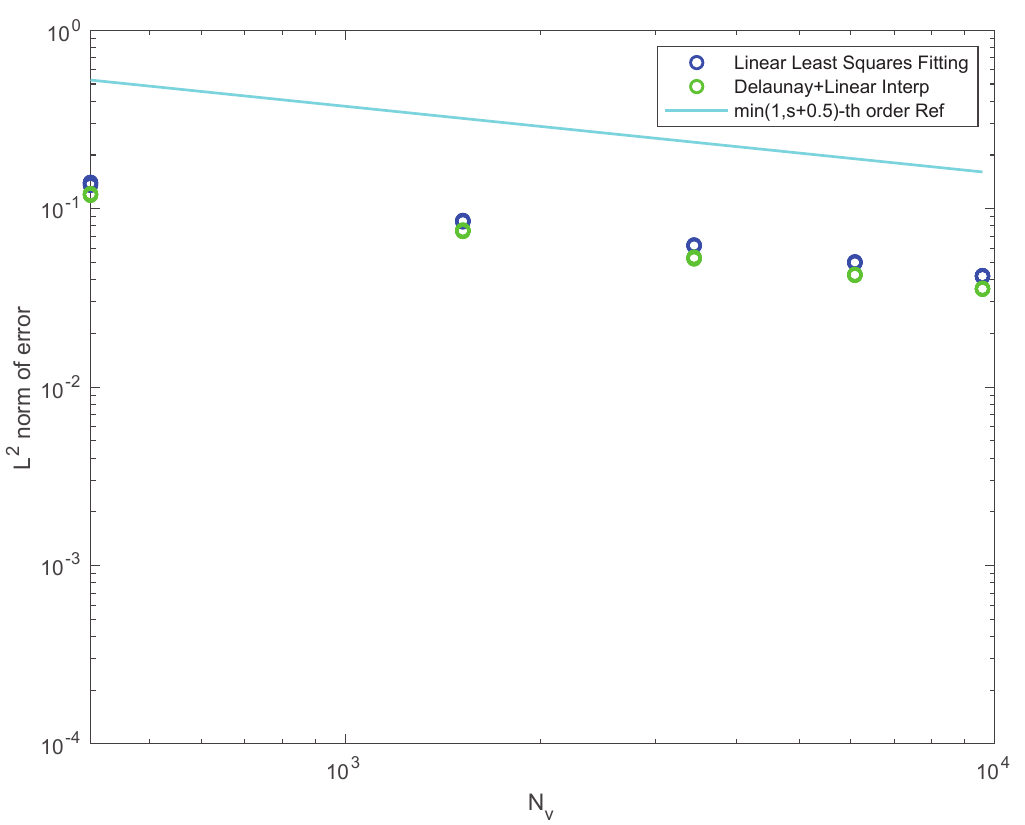}
}
\subfigure[$s=0.5$]{
\includegraphics[width=0.27\linewidth]{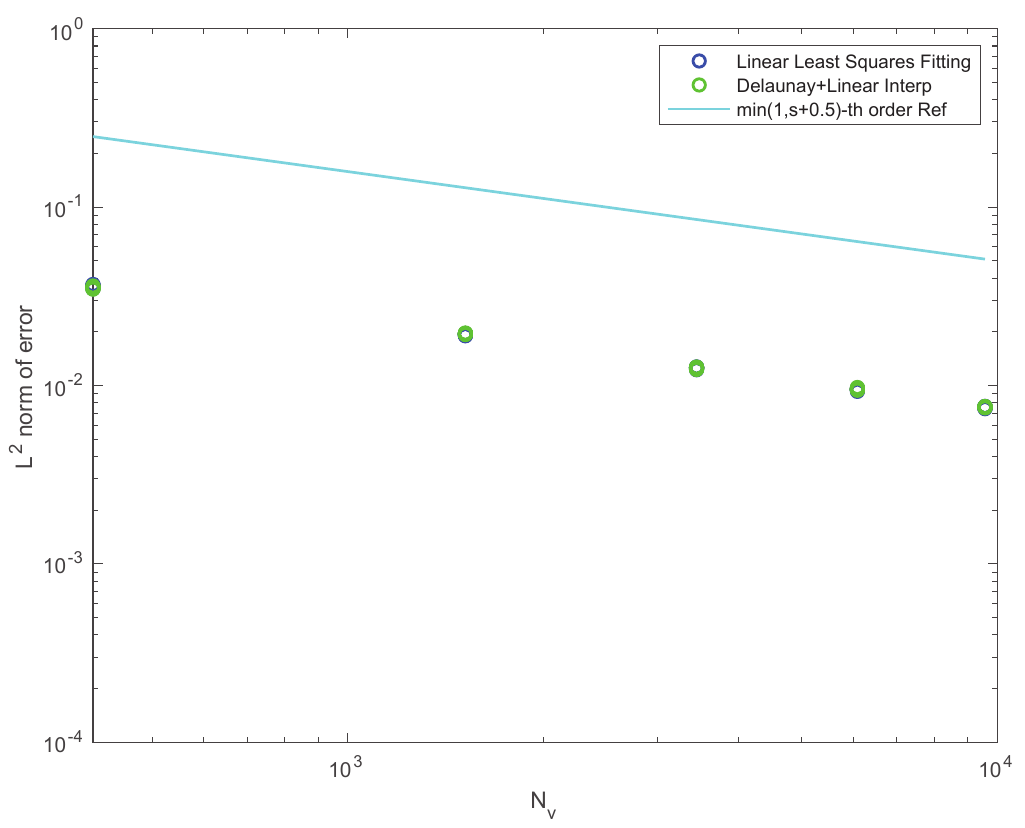}
}
\subfigure[$s=0.75$]{
\includegraphics[width=0.27\linewidth]{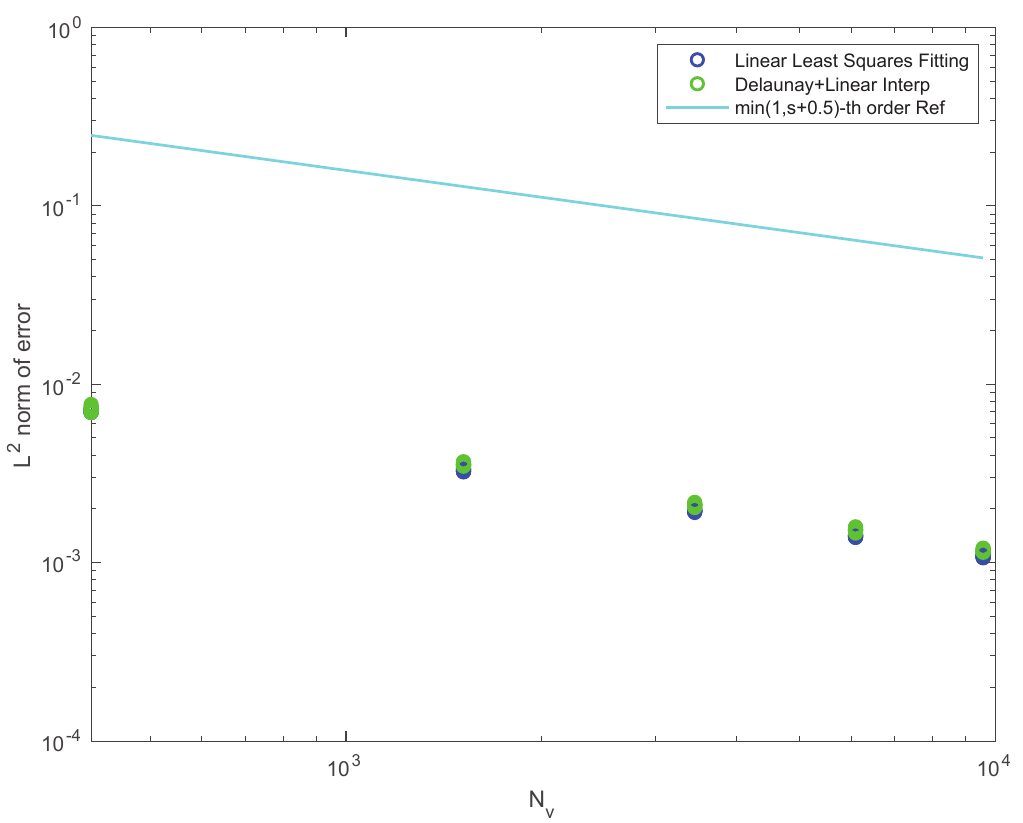}
}
\caption{Example~\ref{exam-1}. The solution error in $L^2$ norm is plotted as a function of $N_v$ for perturbed Point Cloud 2.}
\label{fig:Exam-1-4}
\end{figure}

\begin{figure}[ht!]
\centering
\subfigure[Point Cloud 3]{
\includegraphics[width=0.27\linewidth]{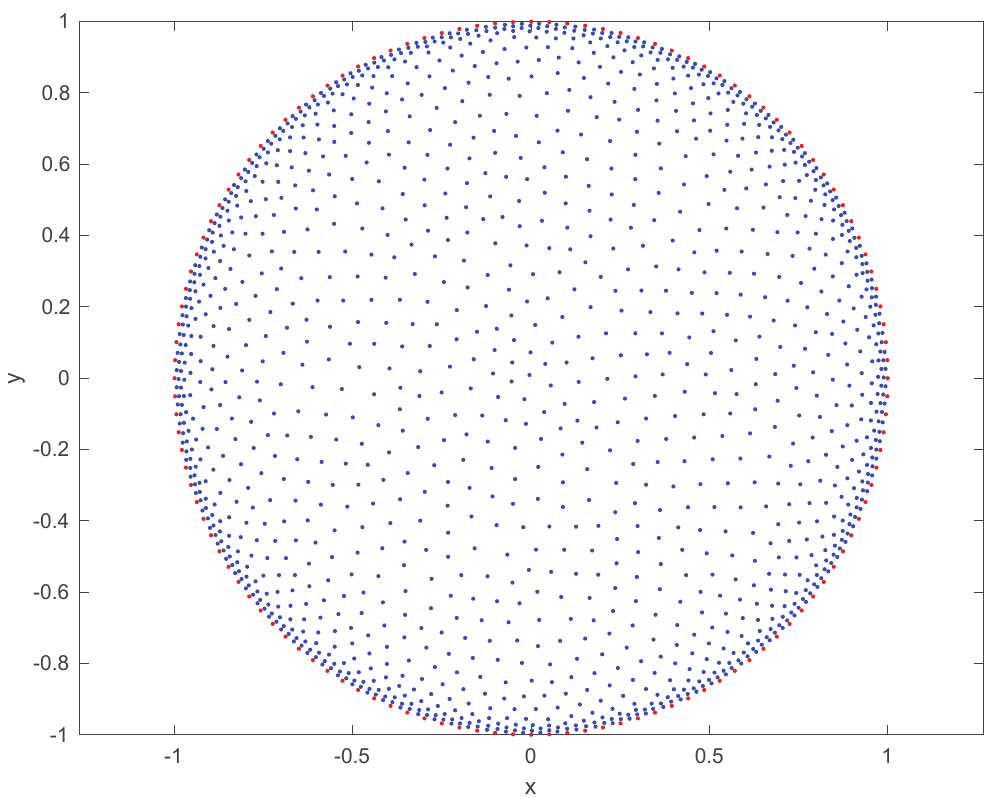}
\label{fig:Exam-1-5-a}
}
\subfigure[Point Cloud 4]{
\includegraphics[width=0.27\linewidth]{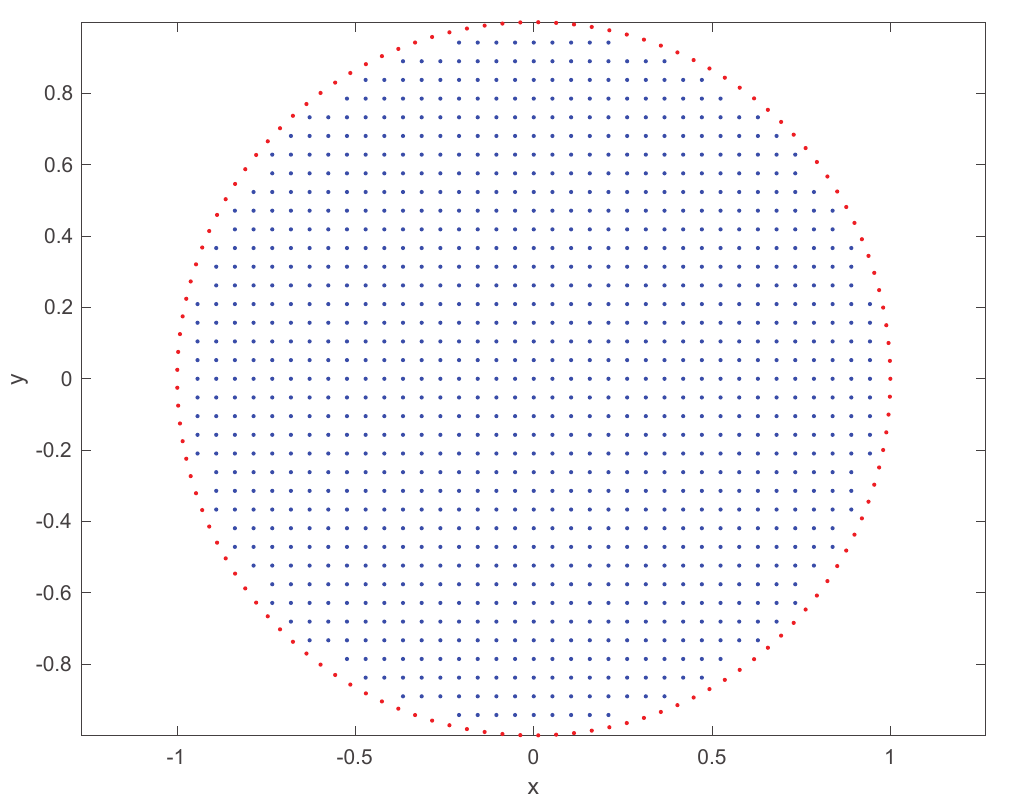}
}
\caption{Example~\ref{exam-1}. Point Clouds 3 and 4.}
\label{fig:Exam-1-5}
\end{figure}

\begin{figure}[ht!]
\centering
\subfigure[$s=0.25$]{
\includegraphics[width=0.305\linewidth]{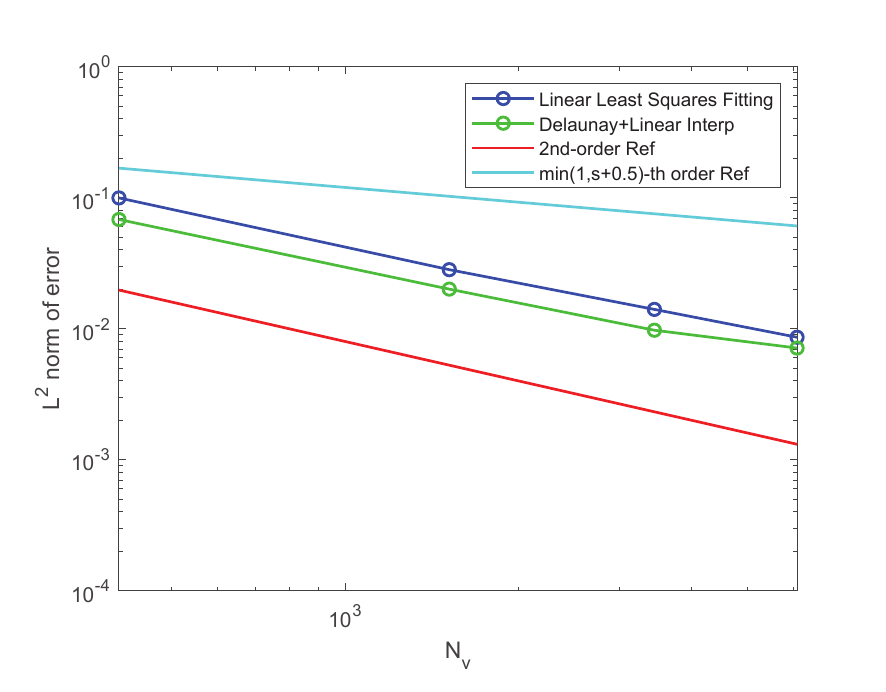}
}
\subfigure[$s=0.5$]{
\includegraphics[width=0.27\linewidth]{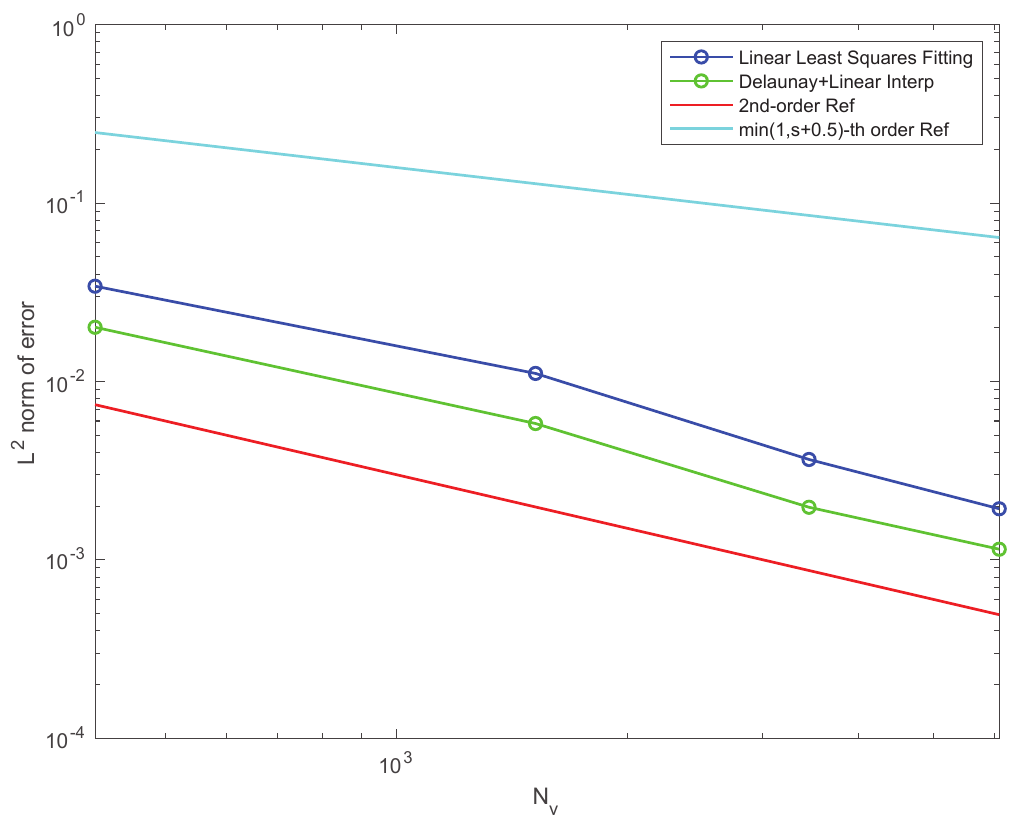}
}
\subfigure[$s=0.75$]{
\includegraphics[width=0.27\linewidth]{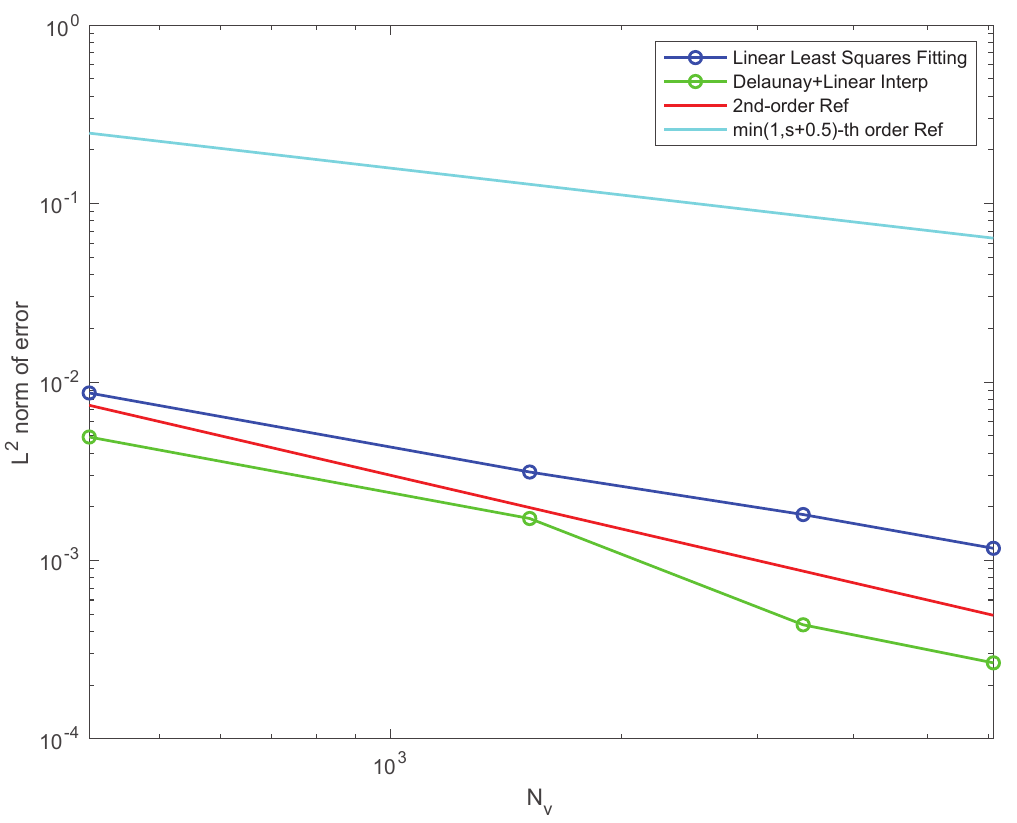}
}
\caption{Example~\ref{exam-1}. The solution error in $L^2$ norm is plotted as a function of $N_v$ for Point Cloud 3.}
\label{fig:Exam-1-6}
\end{figure}

\begin{figure}[ht!]
\centering
\subfigure[$s=0.25$]{
\includegraphics[width=0.27\linewidth]{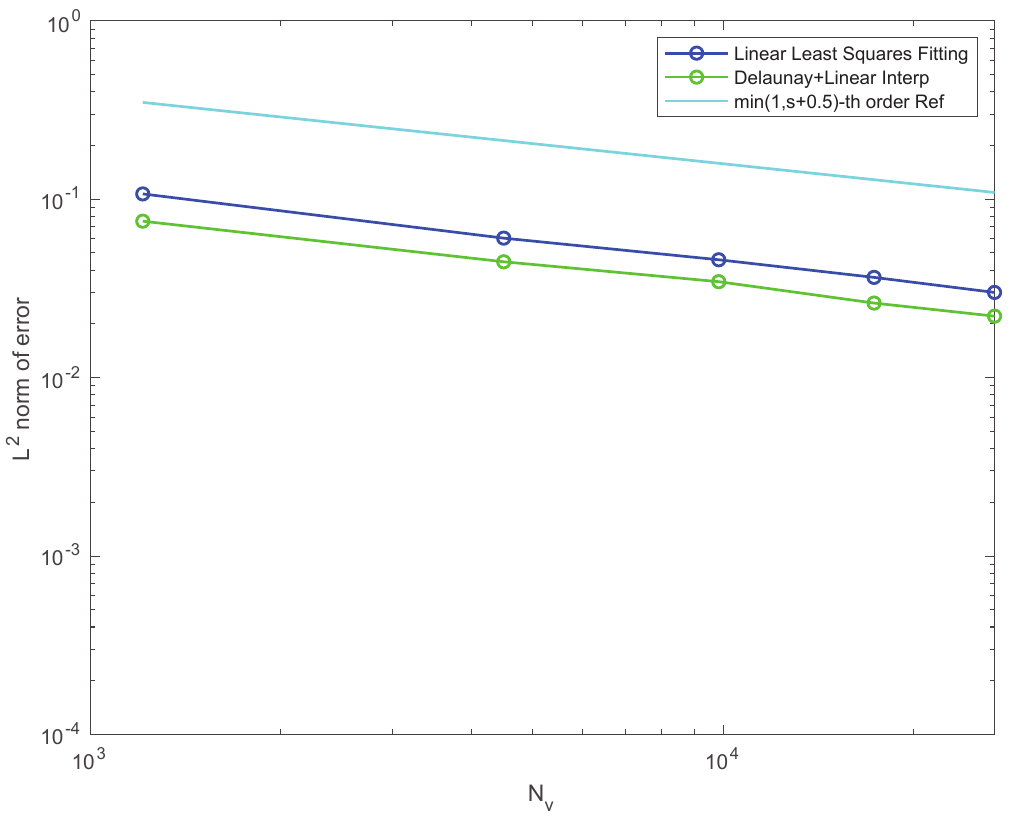}
}
\subfigure[$s=0.5$]{
\includegraphics[width=0.27\linewidth]{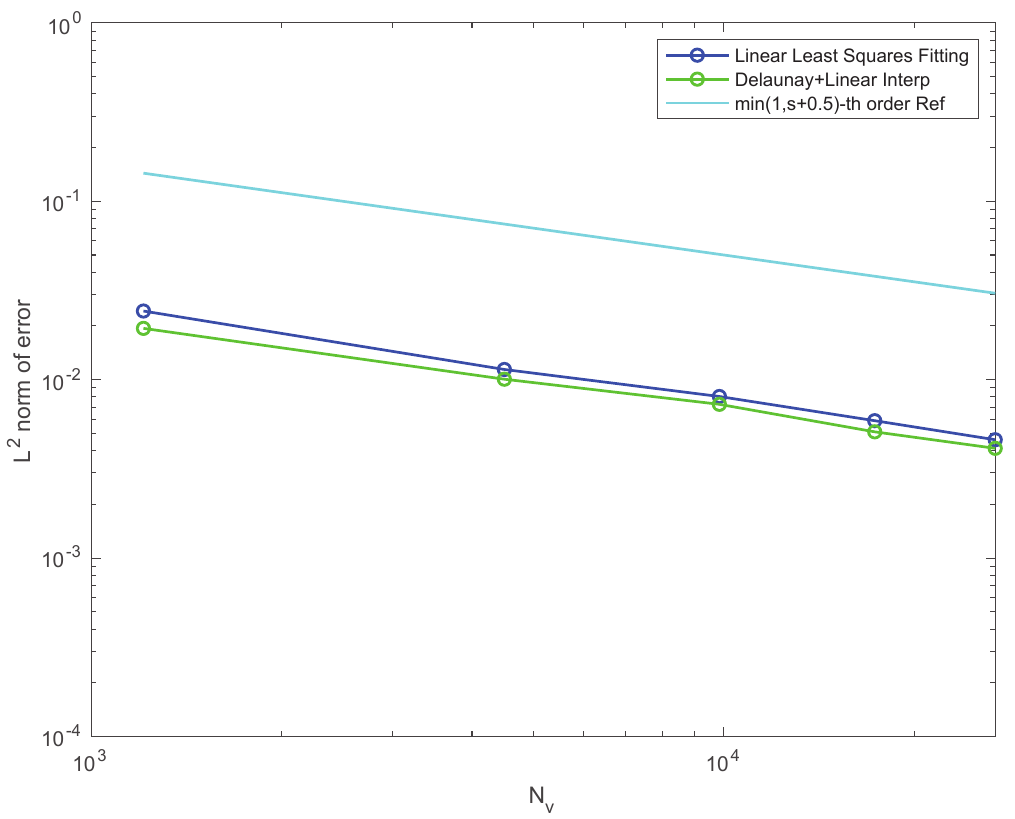}
}
\subfigure[$s=0.75$]{
\includegraphics[width=0.27\linewidth]{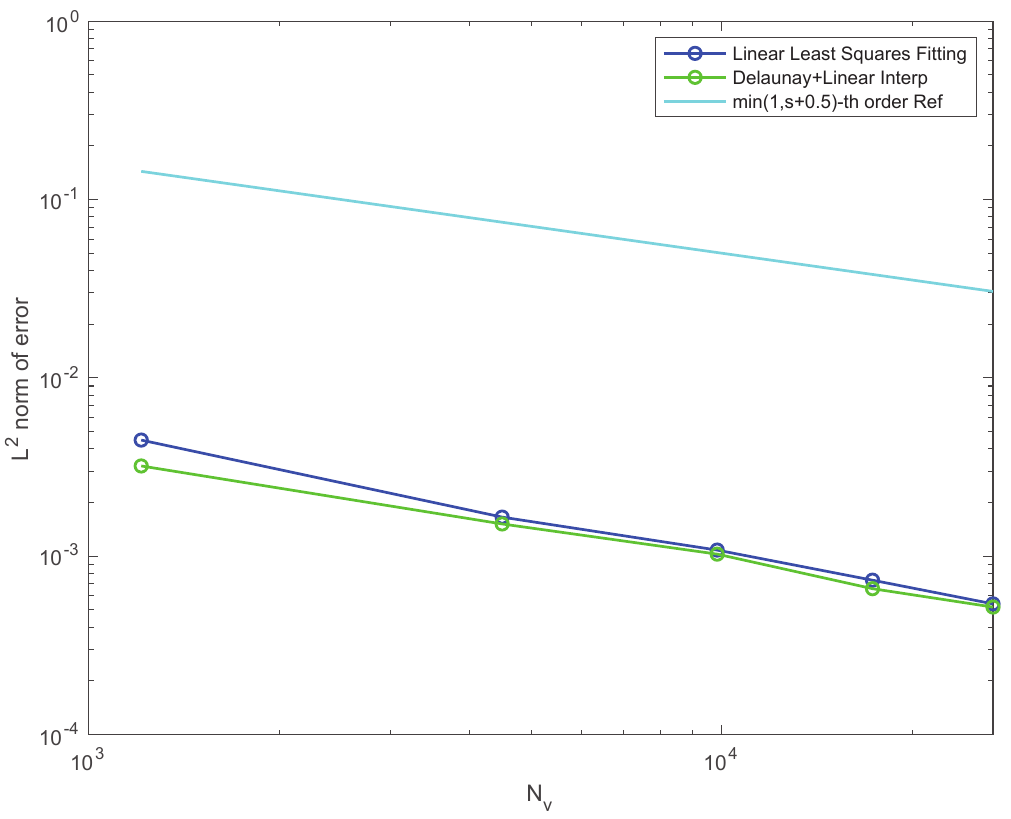}
}
\caption{Example~\ref{exam-1}. The solution error in $L^2$ norm is plotted as a function of $N_v$ for Point Cloud 4.}
\label{fig:Exam-1-7}
\end{figure}

\begin{exam}
\label{exam-1}
In this example, we take $\Omega$ as the unit disk at the origin. In this case,
BVP (\ref{BVP-2}) has the exact solution \cite{Dyda2017}
\begin{equation}
u = \frac{1}{2^{2 s} \Gamma(s+1) \Gamma(s + 1)} (1 - \mid \mbox{\boldmath \( x \)} \mid ^2)_{+}^{s} ,
\label{u-1}
\end{equation}
{ where $(\cdot)_{+}$ stands for taking the positive part.}
The solution is shown in Fig.~\ref{fig:Exam-1-0} for $s = 0.25$, 0.5, and 0.75.

We first consider Point Cloud 1: for any given positive integer $J$,
\begin{equation}
\label{point-set-1}
(x_{j,k}, y_{j,k}) = (\cos(\frac{2 \pi k}{n-1}), \sin(\frac{2 \pi k}{n-1})), \quad k = 0, ..., n-1 ,
\end{equation}
where $j = 0, ..., J-1$ and $n = [j \pi] + 1$.
To examine the robustness of GoFD with respect to point distribution, each interior point in this point cloud is perturbed randomly with noise level
$0.4 \bar{h}$, where $\bar{h} \equiv 1/\sqrt{N_v}$ can be viewed as the average distance among the points.
A point cloud and a corresponding perturbed point cloud are showed in Fig.~\ref{fig:Exam-1-1}.
The solution error in $L^2$ norm is shown in Fig.~\ref{fig:Exam-1-2} as a function of $N_v$.
It can be seen that both the Delaunay and moving least squares fitting approaches lead to comparable error.
Moreover, the convergence behaves like $\mathcal{O}(\bar{h}^{\min(1,s+0.5)})$.
The same convergence rate was observed for GoFD in \cite{Huang2023} for quasi-uniform meshes and
proved theoretically for finite element approximation \cite{Acosta2017}.
More interestingly, the error spreads out only slightly with respect to random perturbations of the point location,
indicating that GoFD is robust with respect to the point distribution.

Point Cloud 2 and its perturbation are shown in Fig.~\ref{fig:Exam-1-3} and the solution error is shown in Fig.~\ref{fig:Exam-1-4}.
This set of points is obtained by removing the connectivity of a triangular mesh generated with mesh generation software.
Once again, the results show that GoFD is robust with respect to the point distribution for both approaches of transfer matrices.

Point Cloud 3 is shown in Fig.~\ref{fig:Exam-1-5}(a) and obtained by removing the connectivity of an adaptive triangular mesh
generated with the so-called MMPDE moving mesh method \cite{HR11,Huang2023} based on the Hessian of the function
$(\text{dist}((x,y),\partial \Omega))^s$ representing the singularity of the BVP solution, where
$\text{dist}((x,y),\partial \Omega))$ denotes the distance from $(x,y)$ to the boundary of $\Omega$.
The solution error is shown in Fig.~\ref{fig:Exam-1-6}, showing second-order convergence for both approaches for constructing
the transfer matrix.

Point Cloud 4 is shown in Fig.~\ref{fig:Exam-1-5}(b) and formed by a selection of points on the boundary and the points of
a uniform grid of spacing $\bar{h}$ that fall in $\Omega$ and are $0.5 \bar{h}$ away from the boundary.
The grid points close to the boundary are not selected because they may be close to some boundary points and their selection
can lead to a very small $h_{\text{FD}}$. (Recall that $h_{\text{FD}}$ is chosen to be half the minimum distance
among the points in the given point cloud.)
The solution error is shown in Fig.~\ref{fig:Exam-1-7} that shows convergence order $\mathcal{O}(\bar{h}^{\min(1,s+0.5)})$
for both approaches for constructing the transfer matrix.
\qed
\end{exam}

\begin{figure}[ht!]
\centering
\subfigure[Domain and point cloud]{
\includegraphics[width=0.27\linewidth]{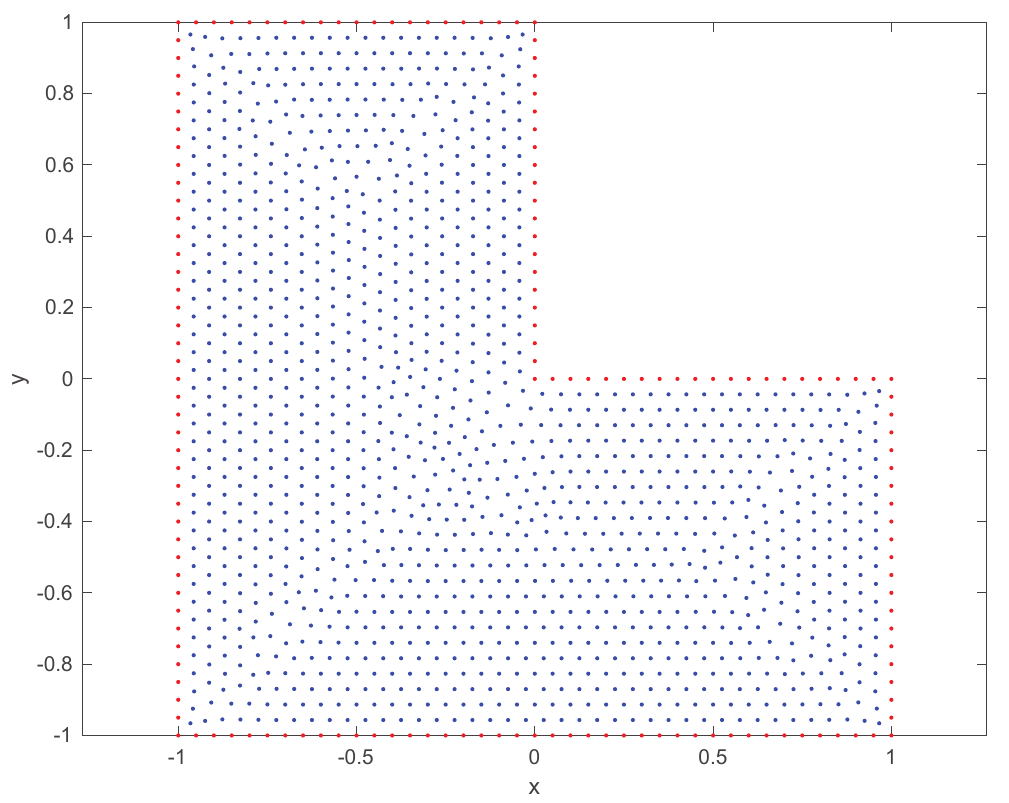}
}
\subfigure[Perturbed point cloud]{
\includegraphics[width=0.31\linewidth]{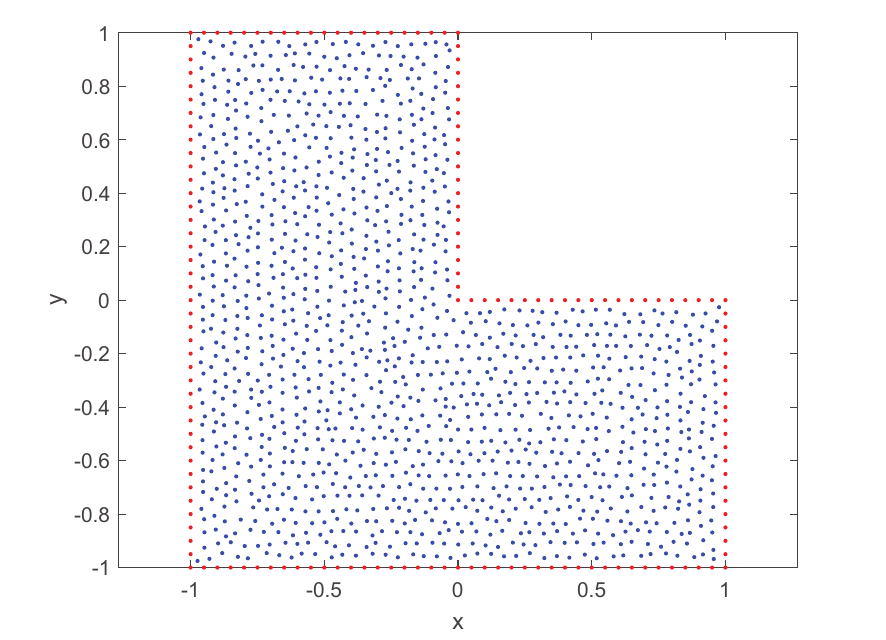}
}
\subfigure[Solution]{
\includegraphics[width=0.27\linewidth]{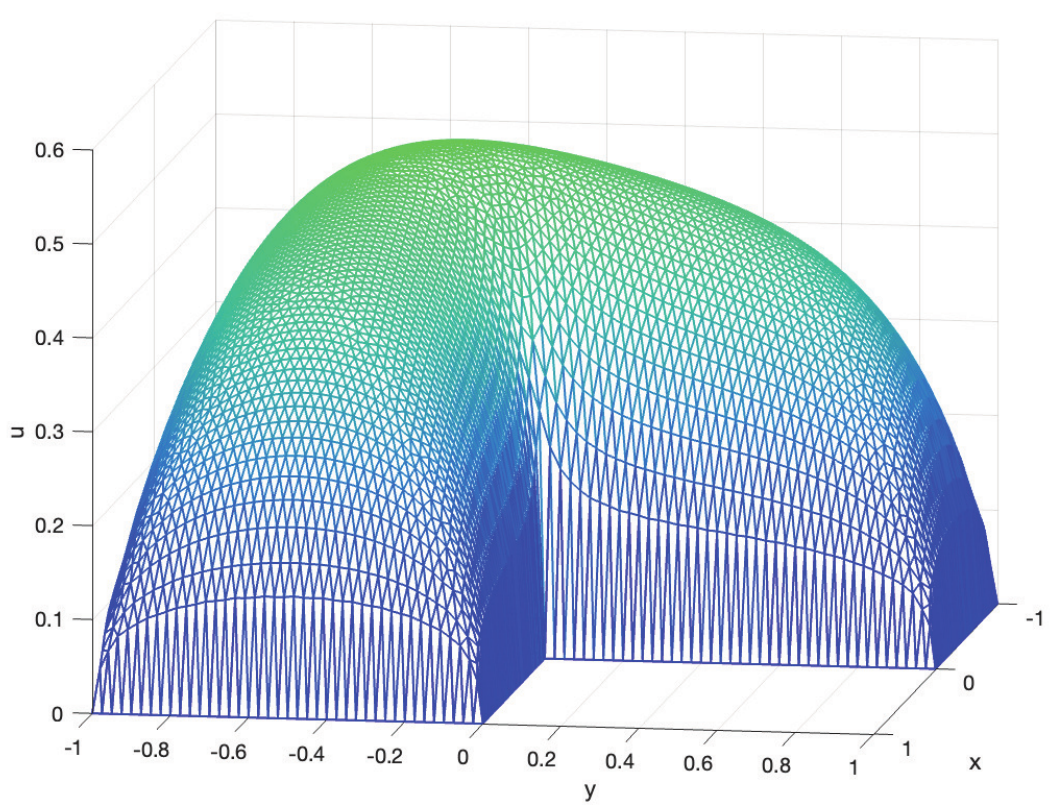}
}
\caption{Example~\ref{exam-2}. Point clouds and solution ($s = 0.5$).}
\label{fig:Exam-2-1}
\end{figure}

\begin{figure}[ht!]
\centering
\subfigure[$s=0.25$]{
\includegraphics[width=0.3\linewidth]{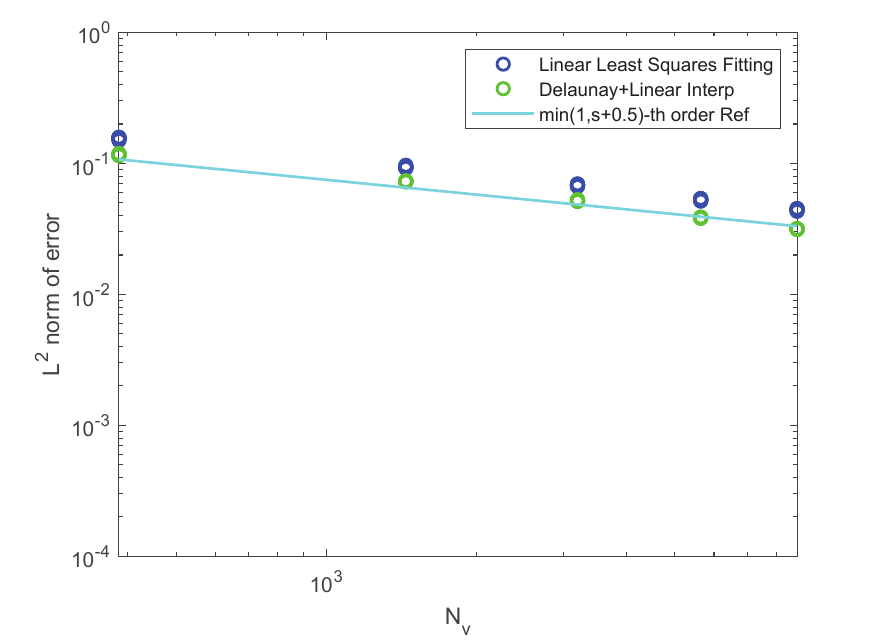}
}
\subfigure[$s=0.5$]{
\includegraphics[width=0.3\linewidth]{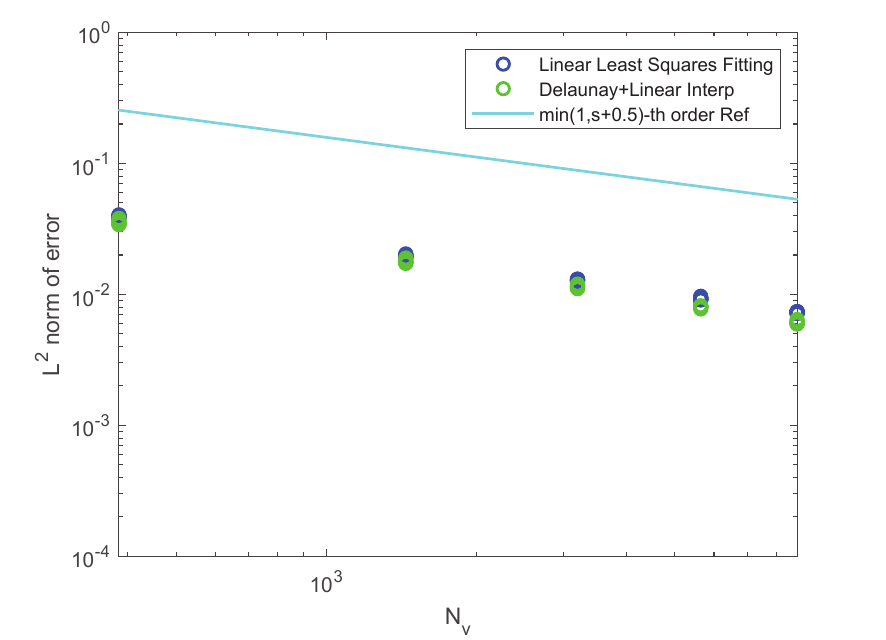}
}
\subfigure[$s=0.75$]{
\includegraphics[width=0.3\linewidth]{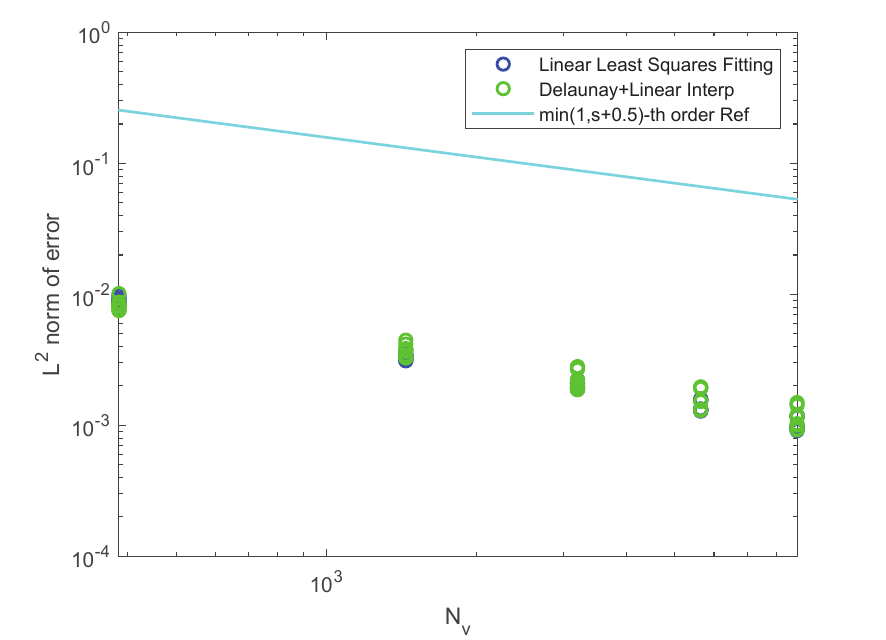}
}
\caption{Example~\ref{exam-2}. The solution error in $L^2$ norm is plotted as a function of $N_v$ for perturbed Point Cloud
shown in Fig.~\ref{fig:Exam-2-1}.}
\label{fig:Exam-2-2}
\end{figure}

\begin{exam}
\label{exam-2}
This example is (\ref{BVP-1}) with $f = 1$ and $\Omega$ being $L$-shaped (cf. BVP (\ref{BVP-2})).
An analytical exact solution is not available for this example. A computed solution obtained with a large number of points
is used as the reference solution in computing the solution error in $L^2$ norm.

A point cloud obtained by removing the connectivity of a unstructured triangular mesh, its perturbation (with noise level $0.4\bar{h}$),
and a computed solution $(s = 0.5)$ are shown in Fig.~\ref{fig:Exam-2-1}.
The solution error in $L^2$ norm is shown as a function of $N_v$ in Fig.~\ref{fig:Exam-2-2}.
The convergence behaves like $\mathcal{O}(\bar{h}^{\min(1,s+0.5)})$ for both the Delaunay and moving least squares fitting approaches,
which is consistent with those of GoFD \cite{Huang2023} and finite element approximation \cite{Acosta2017}
for quasi-uniform meshes.
Moreover, the results show that GoFD is robust with respect to the point distribution for both approaches for constructing the transfer matrix.
\qed
\end{exam}

\begin{figure}[ht!]
\centering
\subfigure[Domain and point cloud]{
\includegraphics[width=0.27\linewidth]{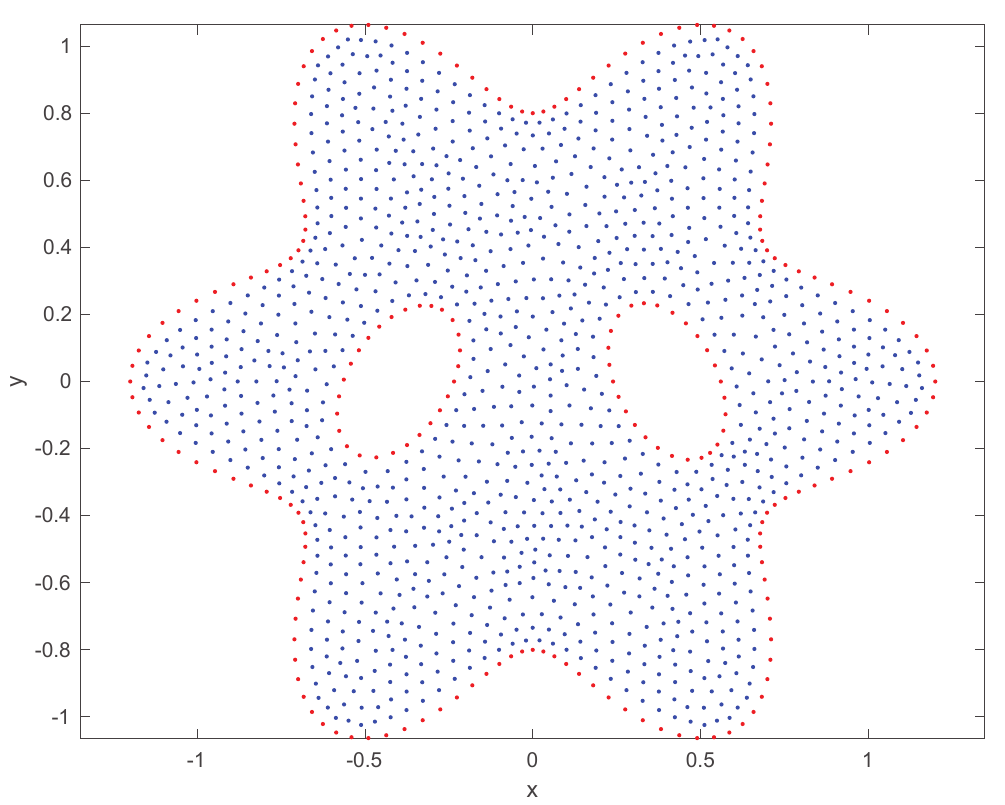}
\label{fig:Exam-3-1-a}
}
\subfigure[Perturbed point cloud]{
\includegraphics[width=0.31\linewidth]{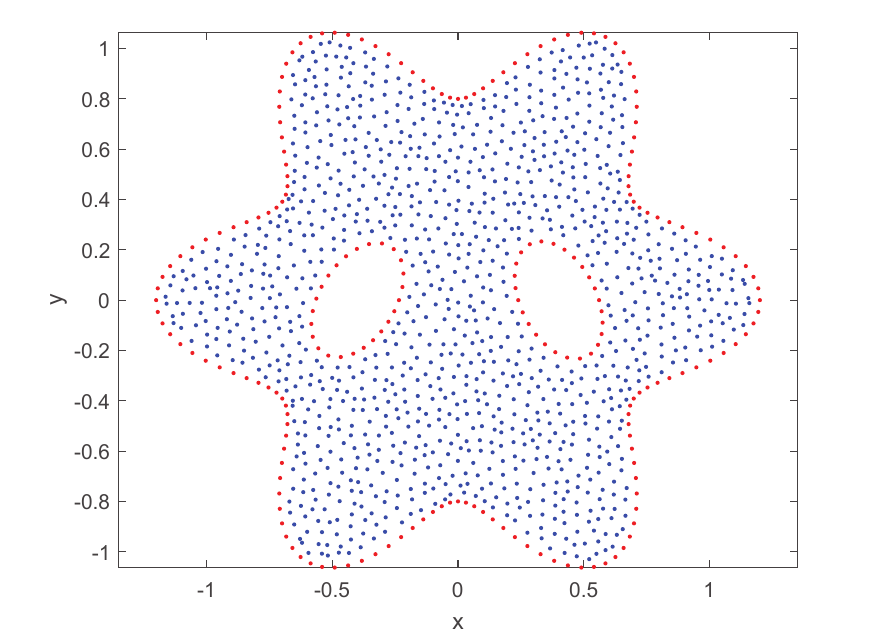}
}
\subfigure[Solution]{
\includegraphics[width=0.27\linewidth]{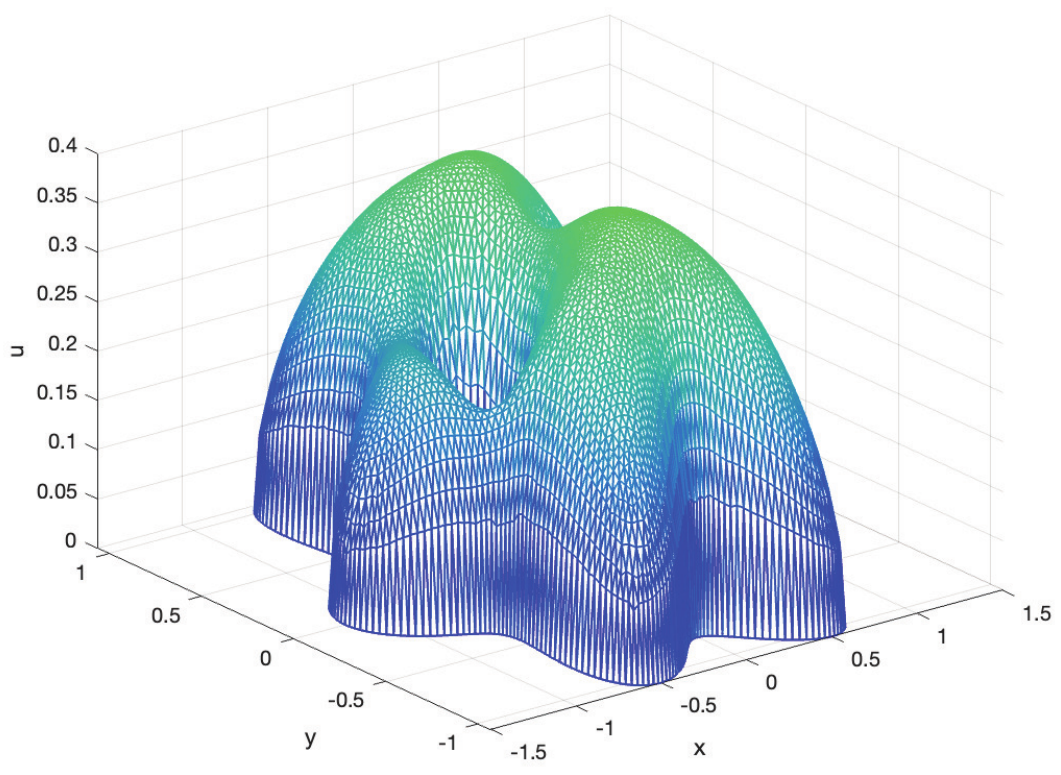}
}
\caption{Example~\ref{exam-3}. Point clouds and solution ($s = 0.5$).}
\label{fig:Exam-3-1}
\end{figure}

\begin{figure}[ht!]
\centering
\subfigure[$s=0.25$]{
\includegraphics[width=0.3\linewidth]{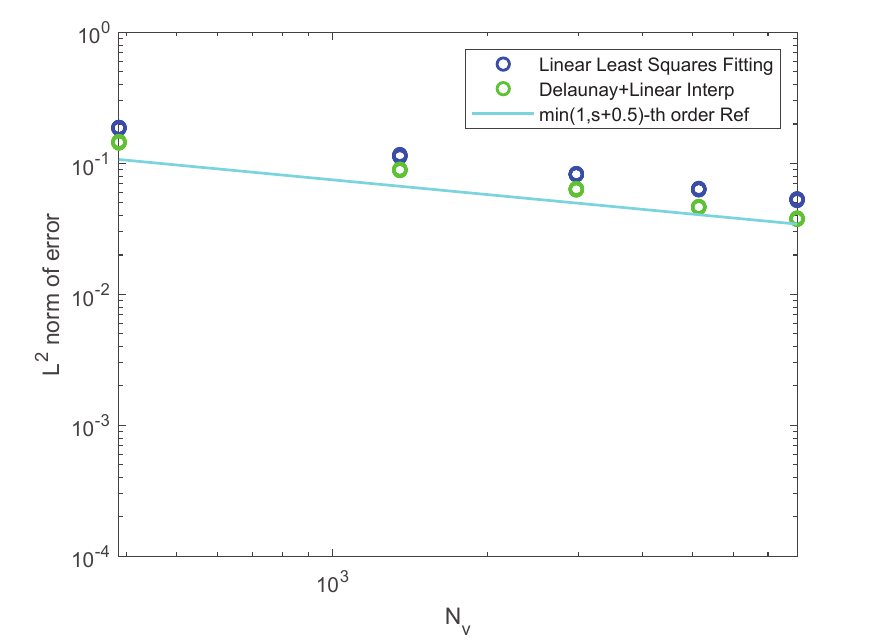}
}
\subfigure[$s=0.5$]{
\includegraphics[width=0.3\linewidth]{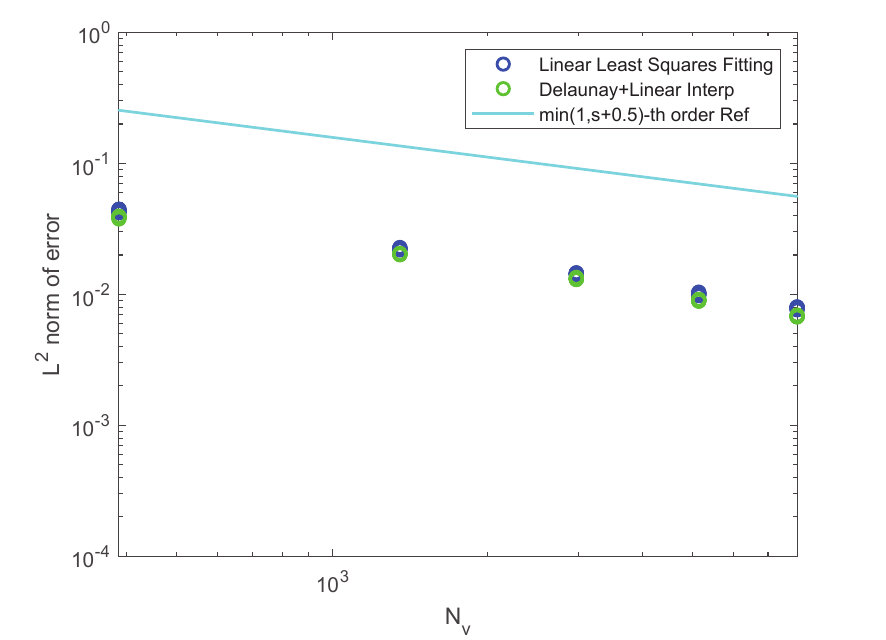}
}
\subfigure[$s=0.75$]{
\includegraphics[width=0.3\linewidth]{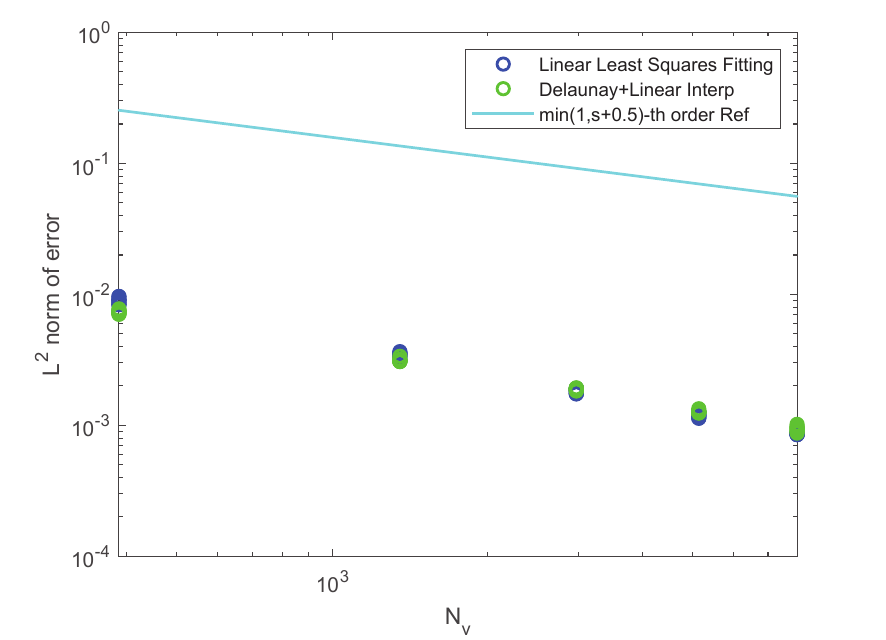}
}
\caption{Example~\ref{exam-3}. The solution error in $L^2$ norm is plotted as a function of $N_v$ for perturbed Point Cloud shown in Fig.~\ref{fig:Exam-3-1}.}
\label{fig:Exam-3-2}
\end{figure}

\begin{figure}[ht!]
\centering
\subfigure[$s=0.25$, Delaunay]{
\includegraphics[width=0.3\linewidth]{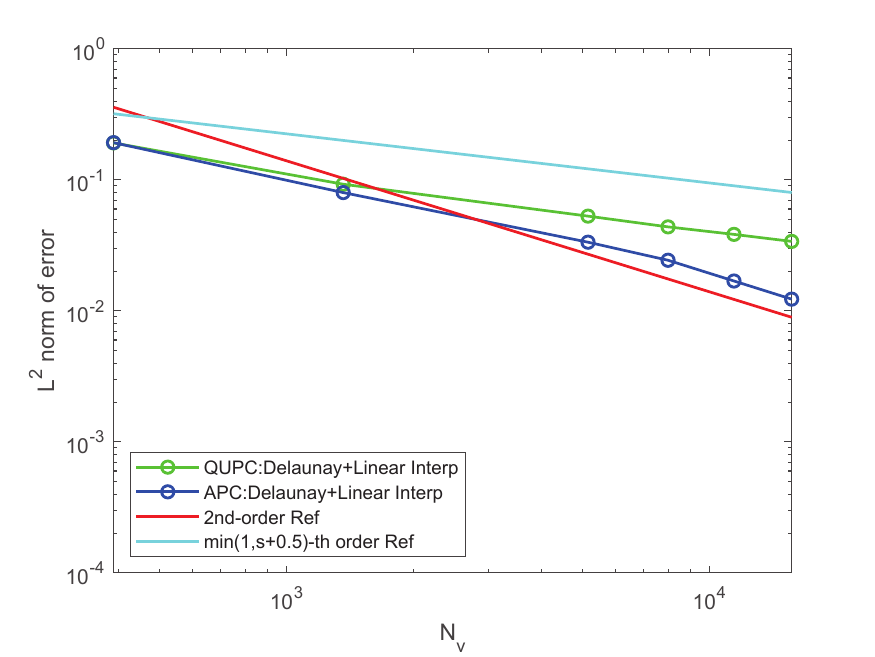}
}
\subfigure[$s=0.5$, Delaunay]{
\includegraphics[width=0.3\linewidth]{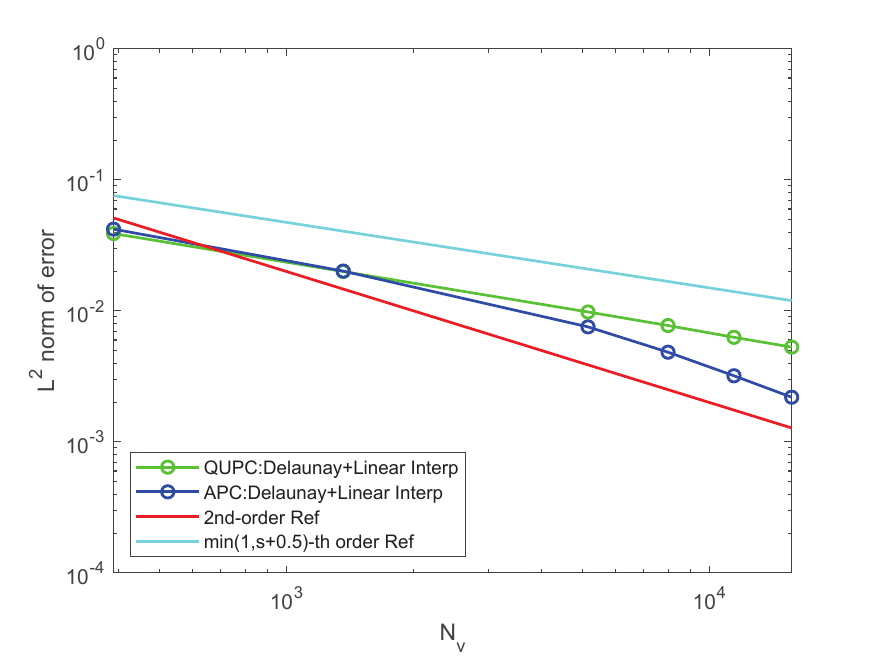}
}
\subfigure[$s=0.75$, Delaunay]{
\includegraphics[width=0.3\linewidth]{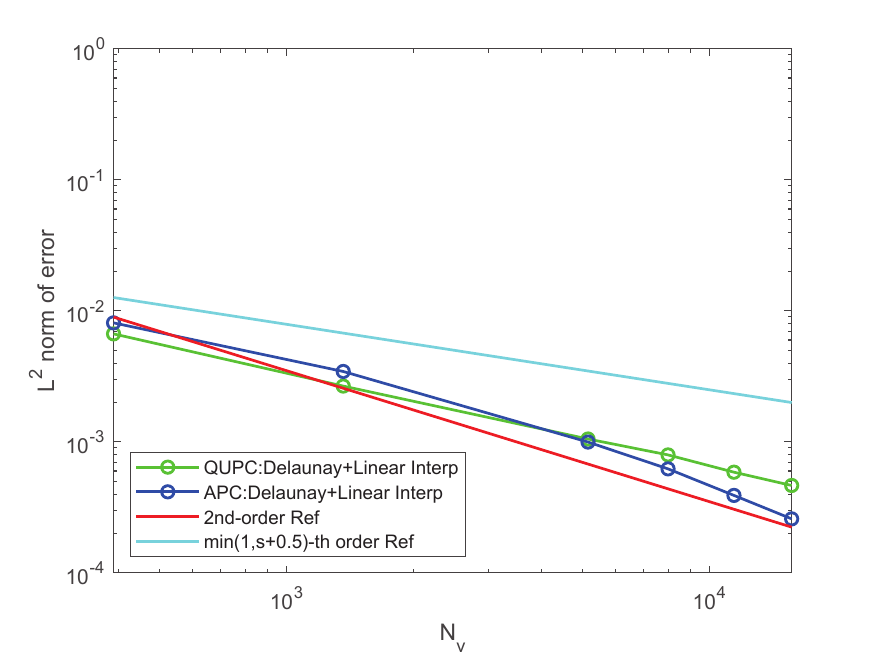}
}
\subfigure[$s=0.25$, Least squares]{
\includegraphics[width=0.3\linewidth]{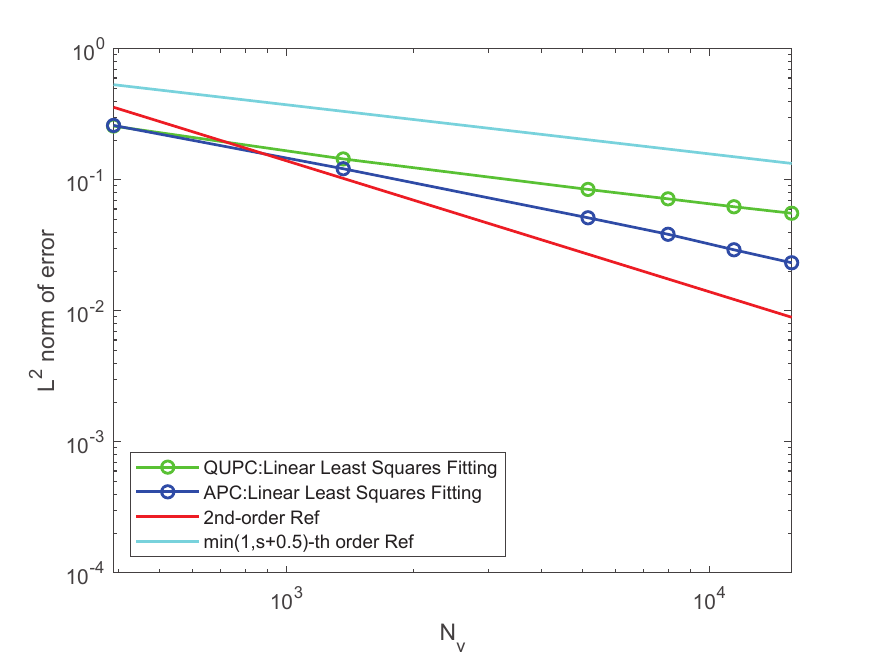}
}
\subfigure[$s=0.5$, Least squares]{
\includegraphics[width=0.3\linewidth]{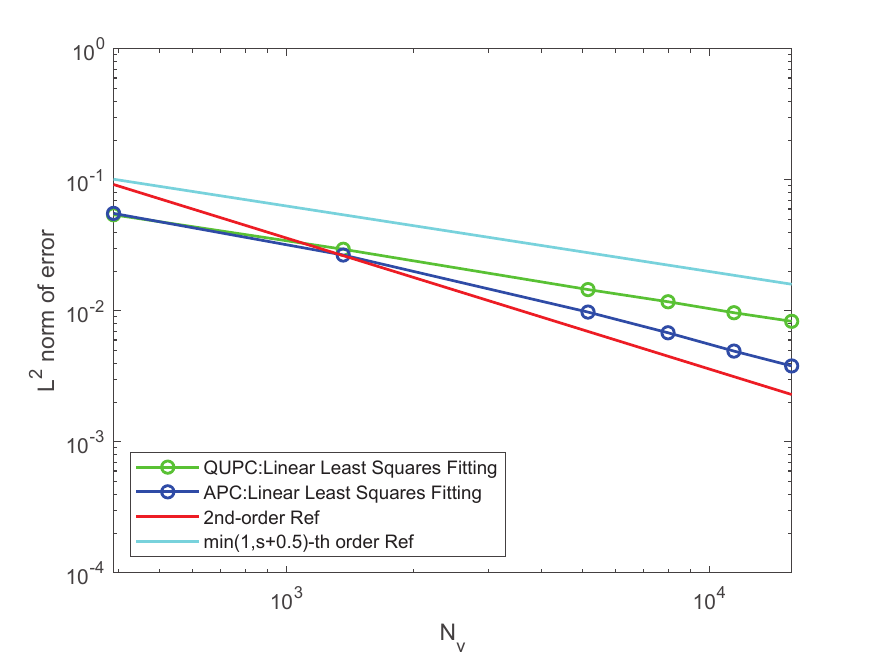}
}
\subfigure[$s=0.75$, Least squares]{
\includegraphics[width=0.3\linewidth]{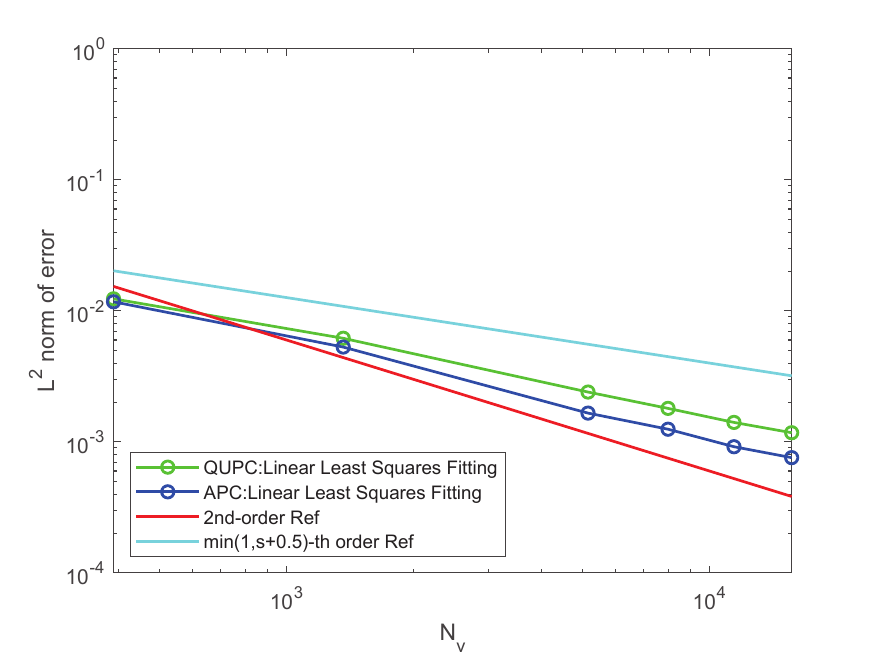}
}
\caption{Example~\ref{exam-3}. The solution error in $L^2$ norm is plotted as a function of $N_v$
for the Delaunay + linear interpolation and linear least squares fitting approaches
for QUPC (quasi-uniform point cloud, cf. point cloud in Fig.~\ref{fig:Exam-3-1-a}) and
APC (adaptive point cloud, cf. point cloud in Fig.~\ref{fig:Exam-1-5-a}).}
\label{fig:Exam-3-3}
\end{figure}

\begin{exam}
\label{exam-3}
This last example is (\ref{BVP-1}) with $f = 1$ and $\Omega$ as shown in Fig.~\ref{fig:gridoverlay-1} (cf. BVP (\ref{BVP-2})).
The geometry of $\Omega$ is complex, with the wavering outer boundary and two holes inside.
An analytical exact solution is not available for this example. A computed solution with a large number of points
is used as the reference solution in computing the $L^2$ norm of the solution error.

A point cloud obtained by removing the connectivity of an unstructured triangular mesh, its perturbation with noise level $0.4 \bar{h}$,
and a computed solution $(s = 0.5)$ are shown in Fig.~\ref{fig:Exam-3-1}.
The solution error is shown in Fig.~\ref{fig:Exam-3-2}. The convergence is about $\mathcal{O}(\bar{h}^{\min(1,s+0.5)})$
for both the Delaunay and least squares fitting approaches of constructing the transfer matrix.
Moreover, the error spreads out only slightly for perturbed point clouds of the same number of points, indicating that GoFD with
the transfer matrix constructed using both approaches is robust with respect to the point distribution.
{
The $L^2$ error is plotted in Fig.~\ref{fig:Exam-3-3} for both approaches and both quasi-uniform and adaptive point clouds.
The results show the second-order convergence for the Delaunay approach and nearly second-order convergence for the linear
least squares fitting approach when adaptive point clouds are used.
}
\qed
\end{exam}

\section{Conclusions and further remarks}
\label{SEC:conclusions}

In the previous sections we have studied the meshfree solution of homogeneous Dirichlet problems
of the fractional Laplacian through the recently developed grid-overlay finite difference method (GoFD) \cite{Huang2023}.
A key to the success of GoFD in a meshfree setting is the construction of the transfer matrix for a given point cloud
for the underlying domain. In this work we have proposed two approaches, one based on the moving least squares fitting
with inverse distance weighting
and the other based on the Delaunay triangulation. While not completely meshfree, the Delaunay triangulation approach
is still a worthwhile option due to its wide use and success in data processing for scattered data interpolation and
low mesh quality requirement by the transfer matrix in GoFD.

Numerical results have been presented for a selection of  two-dimensional examples with convex and concave domains
and various types of point clouds. They have demonstrated that the two approaches for the transfer matrix lead to comparable
solution error. Moreover, as the number of points increases, the solution error in $L^2$ norm
behaves like $\mathcal{O}(\bar{h}^{\min(1,0.5+s)})$ for quasi-uniform point clouds and second-order for adaptive point clouds,
where $\bar{h} = 1/\sqrt{N_v}$ and $N_v$ is the number of points in the cloud.
This is consistent with observations made \cite{Huang2023} for GoFD with unstructured meshes.
Furthermore, the results obtained with a series of random perturbations to the point location show that GoFD with both approaches
for constructing the transfer matrix is robust with respect to the point distribution.

Finally, we remark that GoFD and the approaches studied here for constructing the transfer matrix can be used in other dimensions
without major modifications. Moreover, it has been shown in \cite{Huang2023} that mesh adaptation is necessary to improve
computational accuracy and convergence order for fractional Laplacian boundary value problems.
These improvements have also been observed in Examples~\ref{exam-1} and \ref{exam-3}.
Combining GoFD with meshfree adaptation strategies (e.g., see \cite{Somasekhar-2012,Suchde-2023})
is an interesting topic for future research.

\vspace{20pt}

\noindent
\textbf{Acknowledgments.} W.H. was supported in part by the University of Kansas General Research Fund FY23
and the Simons Foundation through grant MP-TSM-00002397.
J.S. was supported in part by the National Natural Science Foundation of China through grant [12101509].

%

\end{document}